\documentclass[12pt]{article}
\usepackage[pctex32]{graphics}
\usepackage{amsthm,amsmath,amssymb,amscd,verbatim,epsfig}
\usepackage{amssymb}
\usepackage{graphicx}
\newcommand{\beq}{\begin{equation}}
\newcommand{\eeq}{\end{equation}}

\date{}
\newcommand{\al}{\alpha}

\newcommand{\be}{\beta}
\newcommand{\ld}{\lambda}
\newcommand{\si}{\sigma}

\newcommand{\dl}{\delta}

\newcommand{\ra}{\rightarrow}

\newcommand{\sub}{\subset}
\newcommand{\supp}{\supset}

\newcommand{\sq}{$\blacksquare$}

\newcommand\om{\omega}

\newcommand\tha{\theta}

\newcommand\vp{\varphi}
\newcommand\vpm{\varphi_\mu}
\newcommand\vpr{\varphi_\rho}
\newcommand\vps{\varphi_\si}

\newcommand\vpl{\varphi_\ld}

\newcommand\vpo{\varphi_\om}

\begin{document}
\title{On\ the\ one-sided\ and\ two-sided\ similarities\ or\ weak\ similarities\ of\ permutations}
\author{Bau-Sen Du \\ [.5cm]
Institute of Mathematics \\
Academia Sinica \\
Taipei 11529, Taiwan \\
dubs@math.sinica.edu.tw \\}
\maketitle
\begin{abstract}
Let $n \ge 3$ be an integer.  Let $P_n = \{ 1, 2, 3, \cdots, n-1, n \}$ and let $S_n$ be the symmetric group of permutations on $P_n$.  Motivated by the theory of discrete dynamical systems on the interval, we associate each permutation $\si_n$ in $S_n$ a (zero-one) Petrie matrix $M_{\si_n,n-1}$ in $GL(n-1,{{\mathbb{R}}})$ (which is generally not the same as the usual permutation matrix).  Then, for any two permutations $\si_n$ and $\rho_n$ in $S_n$, the notions of right, left and two-sided similarities (and weak similarities respectively) of $\si_n$ and $\rho_n$ are introduced using the similarities (and the characteristic polynomials respectively) of the correspnding Petrie matrices of some extended permutations related to $\si_n$ and $\rho_n$ and examples are presented.  As a by-product, we obtain ways to construct countably infinitely many pairs of Petrie matrices which are similar.  
\end{abstract}

\section{Introduction}
For every positive integer $i$, we shall always let $J_i = [i, i+1]$ and, for every integer $m \ge 2$, we shall always let  $W_{{{\mathbb{F}}}^m} = \biggl\{ \, \sum\limits_{i=1}^m r_iJ_i  :  r_i \in {{\mathbb{F}}}$, $1 \le i \le m  \biggr\}$ be the $m$-dimensional vector space over the field ${{\mathbb{F}}}$ with $B_m = \{  J_i : 1 \le i \le m  \}$  as a (standard) basis.  For convenience, for integers $1 \le j < k \le m+1$, we let $[j, k]$ denote {\it the element} $\sum\limits_{i=j}^{k-1} J_i$ in $W_{{{\mathbb{F}}}^m}$.  Let $n \ge 3$ be a fixed integer and let $P_n$ be the set of all integers in $[1, n]$.  Let $\si$ be a map from $P_n$ into itself such that $\si(i) \ne \si(i+1)$ for all $1 \le i \le n-1$.  The map $\si$ defines a linear transformation, which will always be denoted as $\vps$, from $W_{{{\mathbb{F}}}^{n-1}}$ into itself such that $\vps \biggl(\sum\limits_{i=1}^{n-1} r_iJ_i \biggr) = \sum\limits_{i=1}^{n-1} r_i \vps(J_i)$, where $\vps(J_i) = \sum\limits_{j=\si(i)}^{\si(i+1)-1} J_j$ if $\si(i) < \si(i+1)$ and $\vps(J_i) = \sum\limits_{j=\si(i+1)}^{\si(i)-1} J_j$ if $\si(i+1) < \si(i)$.  Let $M_{\si, n-1} = (a_{i,j})$ be the $(n-1) \times (n-1)$ matrix defined by $a_{i,j} = 1$ for all $\al_i \le j \le \be_i-1$, where $\al_i = \min \{ \si(i), \si(i+1) \}$ and $\be_i = \max \{ \si(i), \si(i+1) \}$, and $a_{ij} = 0$ elsewhere.  Note that, when $\si$ is a permutation, this matrix $M_{\si, n-1}$ need not be the same as the usual permutation matrix.  It is clear that such matrices $M_{\si, n-1}$ have entries either zeros or ones such that the ones in each {\it{row}} occur consecutively.  For our purpose, we define a Petrie matrix to be a square matrix whose entries are either zeros or ones such that the ones in each {\it{row}} occur consecutively.  The determinant of a Petrie matrix is known {\bf{\cite{go}}} and easily seen (by induction) to be either $0$ or $\pm 1$.  So, the matrix $M_{\si, n-1}$ induced by the map $\si$ on $P_n$ is a square Petrie matrix.  We shall call this matrix $M_{\si, n-1}$ the Petrie matrix of $\si$.  It is easy to see that, with respect to the standard basis $B_{n-1} = \{ J_i :1 \le i \le n-1 \}$ for $W_{{{\mathbb{F}}}^{n-1}}$, the Petrie matrix $M_{\si, n-1}$ represents $\vps$ on $W_{{{\mathbb{F}}}^{n-1}}$ in the sense that, for every integer $1 \le i \le n-1$, $\vps(J_i) = \sum\limits_{j=1}^{n-1} a_{ij}J_j$.  In particular, for every element $\sum\limits_{i=1}^{n-1} r_iJ_i$ in $W_{{{\mathbb{F}}}^{n-1}}$, we have $\vps(\sum\limits_{i=1}^{n-1} r_iJ_i) = \sum\limits_{j=1}^{n-1} (\sum\limits_{i=1}^{n-1} r_i a_{ij})J_j$ $=  \sum\limits_{j=1}^{n-1} c_jJ_j$, where $c_j = \sum\limits_{i=1}^{n-1} r_i a_{ij}$.  Or, equivalently, in vector form with respect to the standard basis $B_{n-1} = \{ J_i :1 \le i \le n-1 \}$, we have $(c_1, c_2, \cdots, c_{n-1}) = (r_1, r_2, \cdots, r_{n-1}) \cdot M_{\si, n-1}$.  

\begin{figure}[htb]
\centerline{\epsfig{file=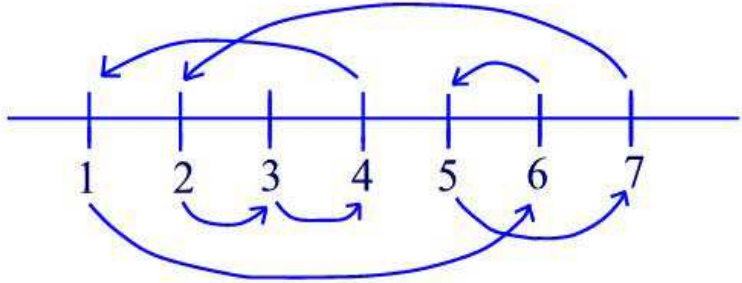,width=12cm,height=5cm}}
\caption{The cyclic permutation $\si_7$.}
\end{figure}

One may ask why anybody would like to study strangely-defined matrices such as the Petrie matrix $M_{\si, n-1}$, especially when conjugate permutations need not have similar corresonding Petrie matrices so defined (cf. {\bf{\cite{ch}}}).  If we look at the Petrie matrix  $M_{\si, n-1}$ from the algebraic point of view, it seems that we have known all interesting things about it {\bf{\cite{aro, gor, go}}}.  So, what else is special about it that warrants the further study?  This answer comes from the dynamical systems point of view.  Although the definition of $M_{\si, n-1}$ may seem opaque, its study {\bf{\cite{al,bl,bl2,du2,du3,ge1,ge2,ge3,ki,ki2,mil,mis}}} is quite common in the theory of discrete dynamical systems on the interval : If we define the linearization of $\si$ on $P_n$ to be the {\it continuous} map $f_\si$ from $[1, n]$ into itself such that $f_\si(k) = \si(k)$ for all integers $1 \le k \le n$ and $f_\si$ is {\it linear} on $J_i$ for all integers $1 \le i \le n-1$ and take $J_i$'s as the vertices of a directed graph and draw an arrow from the vertex $J_i$ to the vertex $J_j$ if $f_{\si}(J_i) \supp J_j$, then $M_{\si, n-1}$ will be the adjacency matrix [$\bf 3$, p.17] of the resulting directed graph.  For example, if $\si_7$ denotes the cyclic permutation $1 \ra 6 \ra 5 \ra 7 \ra 2 \ra 3 \ra 4 \ra 1$ on $P_7$ (see Figure 1), then the adjacency matrix $M_{\si_7, 6}$ of the corresponding directed graph is given as
$$
M_{\si_7, 6} = \left ( \begin{matrix} 
              0 \quad 0 \quad 1 \quad 1 \quad 1 \quad 0 \\
              0 \quad 0 \quad 1 \quad 0 \quad 0 \quad 0 \\
              1 \quad 1 \quad 1 \quad 0 \quad 0 \quad 0 \\
              1 \quad 1 \quad 1 \quad 1 \quad 1 \quad 1 \\
              0 \quad 0 \quad 0 \quad 0 \quad 1 \quad 1 \\ 
              0 \quad 1 \quad 1 \quad 1 \quad 0 \quad 0 \\
              \end{matrix} \right )
$$whose characteristic polynomial is $x^6-3x^5-x^4+5x^3-3x^2-x+1$.  You may have noticed that the coefficients of the  characteristic polynomial of $M_{\si_7, 6}$ are all {\it odd} numbers.  This is no coincidence.  It has been shown in {\bf{\cite{du3}}} that if $\si$ is a {\it cyclic} permutation on $P_n$ then the coefficients of the characteristic polynomial of the adjacency matrix $M_{\si, n-1}$ of the directed graph of $f_\si$ are all {\it odd} numbers and hence nonzero.

This adjacency matrix $M_{\si, n-1}$ turns out to contain many information on the dynamical properties of the map $f_\si$ such as its topological entropy {\bf{\cite{al, bl}}} and Artin-Mazur zeta function {\bf{\cite{ar}}} which are topological invariants.  In {\bf{\cite{du1,du2}}}, when we study the Artin-Mazur zeta function of $f_\si$ for some cyclic permutation $\si$, we are surprised to see that, for some cyclic permutations $\si$ and $\rho$, although $f_\si$ and $f_\rho$ are not topologically conjugate, they still have the same Artin-Mazur zeta function.  Further study reveals that it is because the adajcency matrices of their respective directed graphs are similar.  When we re-examined Theorem 2 of {\bf{\cite{du1}}} a few years ago, we found something interesting hidden behind the theorem (cf. Theorem 7 below).  This initiates the investigation of the one-sided and two-sided similarities and weak similarities of permutations.  Our results on permutations may be complementary from the perspective of discrete dynamical systems on the interval to the book {\it Combinatorics on Permutations} by M. B\'ona {\bf{\cite{bo}}}.  

For convenience, in the sequel, for any set $V_m = \{ v_1, v_2, \cdots, v_m \}$ of $m$ vectors in $W_{{{\mathbb{F}}}^m}$, we shall write $v_i = \sum\limits_{j=1}^m b_{ij}J_j$, $1 \le i \le m$, and call this $m \times m$ matrix $(b_{ij})$ the matrix of $V_m$ (with respect to the standard basis $B_m = \{ J_i :1 \le i \le m \}$) and denote it as $M(V_m|B_m)$.  For any finitely many vectors $w_1, w_2, \cdots, w_k$ in $W_{{{\mathbb{F}}}^m}$, we shall let $< w_1, w_2, \cdots, w_k >$ denote the subspace spanned by $w_1, w_2, \cdots, w_k$.  In Theorem 6, we shall let $F$ denote the field ${\mathbb{Z}}_2 = \{ 0, 1 \}$ of two elements and in other cases, we shall let $F$ denote the field ${\mathbb{R}}$ of real numbers.

\section{Definitions of one-sided and two-sided similarities and weak similarities of permutations}
Let $n \ge 3$ be a fixed integer and let $\si_n$ be a permutation on $P_n$.  In section 1, we associate to each $\si_n$ an $(n-1) \times (n-1)$ matrix $M_{\si_n, n-1}$ called the Petrie matrix of $\si_n$.  It is clear that we can define an equivalence relation on the symmetric group $S_n$ of permutations on $P_n$ by simply using the matrix similarity of $M_{\si_n, n-1}$'s.  However, based on Theorems 5 $\&$ 7 below, we take a different approach.  

In this section, we shall present definitions of one-sided and two-sided similarities or weak similarities of permutations in $S_n$ and some of their properties.  At first look, it may appear somewhat strange to define these similarities in such peculiar ways.  However, they are inspired by the examples in Theorem 7 below which is motivated by the theory of discrete dynamical systems on the interval.  In the theory of discrete dynamical systems on the interval, a period-$n$ orbit of a map determines a (cyclic) permutation $\gamma$ in $S_n$ which, in turn, determines an $(n-1) \times (n-1)$ zero-one Petrie matrix.  By extending the permutation $\gamma$ to the right, left, or, both right and left to a permutation in $S_m$ with $m > n$ in a way suggested by Theorem 7 below, we can compare the characteristic polynomials and discuss the matrix similarity of their corresponding extended Petrie matrices.  It is because of this matrix similarity of the associated Petrie matrices, we call two permutaions in $S_n$ one-sided or two-sided similar.  We shall prove their richness by presenting more examples.

In the sequel, we always let $k$ be a fixed integer which is $\ge 3$.  Let $\si_k$ be a permutation on $P_k = \{ 1,2, \cdots, k-1, k \}$.  For any integer $n \ge 1$, let $\si_{k+n}$ be any permutation on $P_{k+n} = \{ 1,2, \cdots, k+n-1, k+n \}$.  We say that $\si_{k+n}$ is a right extension (see Figure 2) of $\si_k$ if
\begin{itemize}
\item[(1)]
$\si_{k+n}(i) = \si_k(i)$ \, for all integers $1 \le i \le k-1$,

\item[(2)]
$\si_{k+n}(t) = \si_k(k) \le k$ \, for some integer $t$ such that $k+1 \le t \le k+n$,

\item[(3)]
$\si_{k+n}(j) > k$ \, for all integers $j \ne t$ and $k \le j \le k+n$.
\end{itemize}

\begin{figure}[htb]
\centerline{\epsfig{file=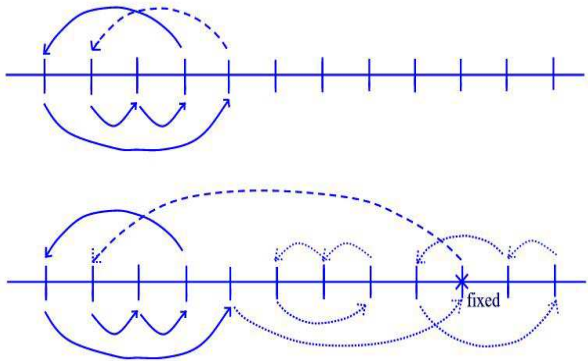,width=8cm,height=5cm}}
\caption{A right extension of the cyclic permutation $\si_5$.}
\end{figure}

\noindent
{\bf Remark.} Let $k$ and $\si_k$ be defined as above.  For any integer $n \ge 1$, let $\beta_{k+1, k+n}$ be any permutation on the set $\{ k+1, k+2, \cdots, k+n \}$.  For any integer $k+1 \le t \le k+n$, we can define a permutation R$(\si_k, \beta_{k+1, k+n}, t)$ (R stands for {\it right}) on the set $\{ 1, 2, 3, \cdots, k+n \}$ by putting (i) R$(\si_k, \beta_{k+1, k+n}, t)(i) = \si_k(i)$ for all $1 \le i \le k-1$; (ii) R$(\si_k, \beta_{k+1, k+n}, t)(k) = \beta_{k+1, k+n}(t)$; (iii) R$(\si_k, \beta_{k+1, k+n}, t)(j) = \beta_{k+1, k+n}(j)$ for all $k+1 \le j \le k+n$ and $j \ne t$; and (iv) R$(\si_k, \beta_{k+1, k+n}, t)(t) = \si_k(k)$.  Then it is easy to see that R$(\si_k, \beta_{k+1, k+n}, t)$ is a right extension of $\si_k$.  Conversely, every right extension of $\si_k$ can be uniquely obtained this way.

Let $\si_k$ and $\rho_k$ be any two permutations on $P_k$ and let $\si_{k+n}$ and $\rho_{k+n}$ be any two permutations on $P_{k+n}$ which are right extensions of $\si_k$ and $\rho_k$ respectively.  We say that the pair $(\si_{k+n}, \rho_{k+n})$ is a synchronized right extension of $\si_k$ and $\rho_k$ if, for some integer $k+1 \le j \le k+n$, $\si_{k+n}(j) = \si_k(k)$ and $\rho_{k+n}(j) = \rho_k(k)$, and $\si_{k+n}(i) = \rho_{k+n}(i)$ for {\it all} integers $i \ne j$ and $k \le i \le k+n$. 

For any two integers $m \ge 1$ and $k \ge 3$, let $\si_k$ and $\si_{m+k}$ be permutations on $P_k$ and on $P_{m+k}$ respectively.  We say that $\si_{m+k}$ is a left extension (see Figure 3) of $\si_k$ if
\begin{itemize}
\item[(1)]
$\si_{m+k}(m+i) = m + \si_k(i)$ \, for all integers $2 \le i \le k$,

\item[(2)]
$\si_{m+k}(s) = m + \si_k(1) \ge m+1$ \, for some integer $s$ such that $1 \le s \le m$,

\item[(3)]
$\si_{m+k}(j) < m+1$ \, for all integers $j \ne s$ and $1 \le j \le m+1$.
\end{itemize}

\noindent
{\bf Remarks.} (1) For simplicity, we sometimes call $\si_{m+k}$ a left extension of the permutation $\si_{m+1,m+k}(x) = \si_k(x-m) + m$ defined on the set $\{ m+1, m+2, m+3, \cdots, m+k-1, m+k \}$ which can be seen as a translation of $\si_k$ to the right by $m$ units.

(2) Let $k$ and $\si_k$ be defined as above.  For any integer $m \ge 1$, let $\alpha_m$ be any permutation on the set $\{ 1, 2, \cdots, m \}$.  For any integer $1 \le s \le m$, we can define a permutation L$(\alpha_m, \si_k, s)$ (L stands for {\it left}) on the set $\{ 1, 2, 3, \cdots, m+k \}$ by putting (a) L$(\alpha_m, \si_k, s)(i) = \alpha_m(i)$ for all $1 \le i \le m$ and $i \ne s$; (b) L$(\alpha_m, \si_k, s)(s) = m+\si_k(1)$; (c) L$(\alpha_m, \si_k, s)(m+1) = \alpha_m(s)$; and (d) L$(\alpha_m, \si_k, s)(j) = m+\si_k(j-m)$ for all $m+2 \le j \le k+m$.  Then it is easy to see that L$(\alpha_m, \si_k, s)$ is a left extension of $\si_k$.  Conversely, every left extension of $\si_k$ can be uniquely obtained this way.

\begin{figure}[htb]
\centerline{\epsfig{file=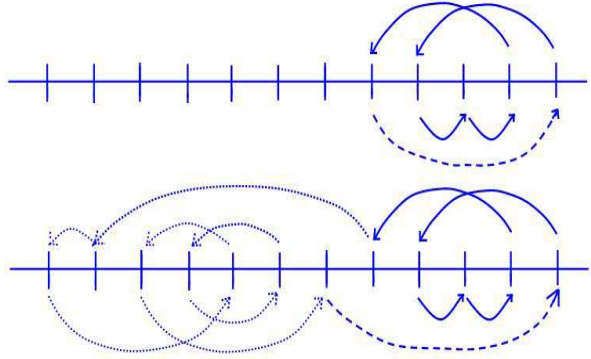,width=8cm,height=5cm}}
\caption{A left extension of the cyclic permutation $\si_5$.}
\end{figure}

\bigskip

Let $\si_k$ and $\rho_k$ be any two permutations on $P_k$ and let $\si_{m+k}$ and $\rho_{m+k}$ be any two permutations on $P_{m+k}$ which are left extensions of $\si_k$ and $\rho_k$ respectively.  We say that the pair $(\si_{m+k}, \rho_{m+k})$ is a synchronized left extension of $\si_k$ and $\rho_k$ if, for some integer $1 \le j \le m$, $\si_{m+k}(j) = m + \si_k(1)$ and $\rho_{m+k}(j) = m + \rho_k(1)$, and $\si_{m+k}(i) = \rho_{m+k}(i)$ for {\it all} integers $i \ne j$ and $1 \le i \le m+1$.

\bigskip

For any integers $m \ge 1$, $k \ge 3$, and $n \ge 1$, let $\si_k$ and $\si_{m+k+n}$ be permutation on $P_k$ and on $P_{m+k+n}$ respectively.  We say that $\si_{m+k+n}$ is a two-sided extension (see Figure 4) of $\si_k$ if
\begin{itemize}
\item[(1)]
$\si_{m+k+n}(m+i) = m + \si_k(i)$ \, for all integers $2 \le i \le k-1$,

\item[(2)]
$\si_{m+k+n}(s) = m + \si_k(1)$ \, for some integer  $s$  such that $1 \le s \le m$,

\item[(3)]
$\si_{m+k+n}(i) < m+1$ \, for all integers $i \ne s$ and $1 \le i \le m+1$,

\item[(4)]
$\si_{m+k+n}(t) = m + \si_k(k)$ \, for some integer  $t$  such that $m+k+1 \le t \le m+k+n$,

\item[(5)]
$\si_{m+k+n}(j) > m+k$ \, for all integers $j \ne t$ and $m+k \le j \le m+k+n$.
\end{itemize}

\begin{figure}[htb]
\centerline{\epsfig{file=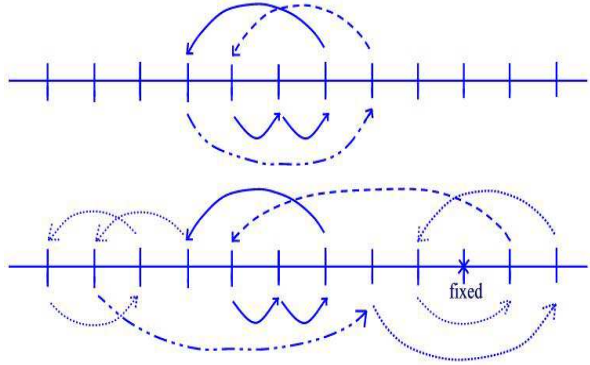,width=8cm,height=5cm}}
\caption{A two-sided extension of the cyclic permutation $\si_5$.}
\end{figure}

\noindent
{\bf Remark.} Let $k$ and $\si_k$ be defined as above.  For any integers $m \ge 1$ and $n \ge 1$, let $\alpha_m$ and $\beta_{m+k+1, m+k+n}$ be any permutations on $P_m$ and on $\{ m+k+1, m+k+2, m+k+3, \cdots, m+k+n \}$ respectively.  For any integers $1 \le s \le m$ and $m+k+1 \le t \le m+k+n$, we can define a permutation T$(\alpha_m, \si_k, \beta_{m+k+1, m+k+n}, s, t)$ (T stands for {\it two-sided}) on the set $P_{m+k+n}$ by putting (a) T$(\alpha_m, \si_k, \beta_{m+k+1, m+k+n}, s, t)(i) = \alpha_m(i)$ for all $1 \le i \le m$ and $i \ne s$; (b) T$(\alpha_m, \si_k, \beta_{m+k+1, m+k+n}, s, t)(s) = m + \si_k(1)$; (c) T$(\alpha_m, \si_k, \beta_{m+k+1, m+k+n}, s, t)(m+1) = \alpha_m(s)$; (d) T$(\alpha_m, \si_k, \beta_{m+k+1, m+k+n}, s, t)(j) = m+\si_k(j-m)$ for all $m+2 \le j \le m+k-1$; (e) T$(\alpha_m, \si_k, \beta_{m+k+1, m+k+n}, s, t)(m+k) = \beta_{m+k+1, m+k+n}(t)$; (f) T$(\alpha_m, \si_k, \beta_{m+k+1, m+k+n}, s, t)(t) = m+\si_k(k)$; and (g) T$(\alpha_m, \si_k, \beta_{m+k+1, m+k+n}, s, t)(j)$ $= \beta_{m+k+1, m+k+n}(j)$ for all $m+k+1 \le j \le m+k+n$ and $j \ne t$.  Then it is easy to see that T$(\alpha_m, \si_k, \beta_{m+k+1, m+k+n}, s, t)$ is a two-sided extension of $\si_k$.  Conversely, every two-sided extension of $\si_k$ can be uniquely obtained this way.

Let $\si_k$ and $\rho_k$ be any two permutations on $P_k$ and let $\si_{m+k+n}$ and $\rho_{m+k+n}$ be any two permutations on $P_{m+k+n}$ which are two-sided extensions of $\si_k$ and $\rho_k$ respectively.  We say that the pair $(\si_{m+k+n}, \rho_{m+k+n})$ is a synchronized two-sided extension of $\si_k$ and $\rho_k$ if, for some integers $1 \le s \le m$ and $m+k+1 \le t \le m+k+n$, we have $\si_{m+k+n}(s) = m + \si_k(1)$, $\rho_{m+k+n}(s) = m + \rho_k(1)$, $\si_{m+k+n}(t) = m + \si_k(k)$, $\rho_{m+k+n}(t) = m + \rho_k(k)$, and $\si_{m+k+n}(i) = \rho_{m+k+n}(i)$ for all integers $i \in [1, m+1] \cup [m+k, m+k+n] \setminus \{ \, s, t \, \}$. 

\bigskip
We now define right, left and two-sided similarities or weakly similarities of any two permutations on $P_k$.  Note that  these definitions have nothing to do with the conjugation in the corresponding symmetric group $S_k$ of permutations on $P_k$.

\noindent
{\bf Definition 1.}
Let $\si_k$ and $\rho_k$ be any two permutations on $P_k$.  We say that $\si_k$ and $\rho_k$ (or the pair $(\si_k, \rho_k)$) are right (left, two-sided respectively) similar if, for {\it every} synchronized right (left, two-sided respectively) extension $(\si, \rho)$ of $\si_k$ and $\rho_k$, the Petrie matrices of $\si$ and $\rho$ are similar.

\noindent
{\bf Definition 2.}
Let $\si_k$ and $\rho_k$ be any two permutations on $P_k$.  We say that $\si_k$ and $\rho_k$ (or the pair $(\si_k, \rho_k)$) are right (left, two-sided respectively) weakly similar if, for {\it every} synchronized right (left, two-sided respectively) extension $(\si, \rho)$ of $\si_k$ and $\rho_k$, the Petrie matrices of $\si$ and $\rho$ have the same characteristic polynomial.

\noindent
{\bf Remarks.} (1) It is easy to see that right (left or two-sided respectively) similarity defines an equivalence relation on the set of permutations on $P_k$ (or, on the symmetric group $S_k$).  

(2) It is easy to see that right (left or two-sided respectively) weakly similarity defines an equivalence relation on the set of permutations on $P_k$ (or, on the symmetric group $S_k$).

(3) There are permutations on $P_k$, for any $k \ge 4$, which are right (left, or two-sided respectively) weakly similar but not similar.  For example, if $\si_k = (34)$ and $\rho_k = (12)(34)$, then $\si_k$ and $\rho_k$ are right weakly similar but not right similar.  That they are right weakly similar can be seen by expanding the determinants $\det(M_{\si_{k+n}, k+n-1} - \ld I)$ and $\det(M_{\rho_{k+n}, k+n-1} - \ld I)$ along the first row using Laplace's formula, where $I$ is the $(k+n-1) \times (k+n-1)$ identity matrix and $(\si_{k+n}, \rho_{k+n})$ is, for any $n \ge 1$, any synchronized right extension of $\si_k$ and $\rho_k$ on $P_{k+n}$.  However, for each $k \ge 4$, $\si_k$ and $\rho_k$ are not right similar.  More examples can be found in section 5.  

In this paper, we shall concentrate only on one-sided or two-sided similarity.  The following results which can be verified easily demonstrate how to obtain more examples of one-sided or two-sided similar permutations from the known ones.

\noindent
{\bf Theorem 1.}
{\it Let $\si_k$ and $\rho_k$ be any two permutations on $P_k$ which are two-sided similar.  Then the following hold:
\begin{itemize}
\item[\rm{(1)}]
For every synchronized two-sided extension $(\si, \rho)$ of $\si_k$ and $\rho_k$, $\si$ and $\rho$ are right, left, and two-sided similar.

\item[\rm{(2)}]
For every synchronized right (left respectively) extension $(\si, \rho)$ of $\si_k$ and $\rho_k$, $\si$ and $\rho$ are left and two-sided (right and two-sided respectively) similar. 
\end{itemize}}

\noindent
{\bf Theorem 2.}
{\it Let $\si_k$ and $\rho_k$ be any two permutations on $P_k$ which are right (left respectively) similar.  Then, for every synchronized right (left respectively) extension $(\si, \rho)$ of $\si_k$ and $\rho_k$, $\si$ and $\rho$ are right (left respectively) similar.}

For any permutation $\theta_k$ on $P_k$, we let $\theta_k^*$ denote the dual permutation of $\theta_k$ on $P_k$ such that $\theta_k^*(i) = k+1 - \theta_k(k+1-i)$, $1 \le i \le k$.  It is easy to see that the dynamics of $\theta_k$ and $\theta_k^*$ on $P_k$ are mirror symmetric.    

\noindent
{\bf Theorem 3.}
{\it Let $\si_k$ and $\rho_k$ be any two permutations on $P_k$ which are right (left, two-sided respectively) similar.  Then $\si_k^*$ and $\rho_k^*$ are left (right, two-sided respectively) similar.}

We can also combine two pairs of permutations with various similarity properties to obtain more examples.  

\noindent
{\bf Theorem 4.}
{\it Let $\ell \ge 4$ and $k \ge \ell+1$ be integers.  Let $\si_\ell$ and $\rho_\ell$ be permutations on $P_\ell$ such that $\si_\ell(\ell) = \ell-1$ and $\si_\ell(\ell-1) < \ell-1$.  Let $\xi$ and $\eta$ be permutations on the set $\{ \ell-1, \ell, \ell+1, \ell+2, \cdots, k-1, k \}$ such that $\xi(\ell) = \ell-1$, and $\eta(\ell-1) = \ell$ and $\eta(\ell) > \ell$.  Let $(\si\xi)_k$, $(\si\eta)_k$, and $(\rho\eta)_k$ be the permutations on $P_k$ defined by 
$$
(\si\xi)_k(x) =   \begin{cases}
               \si_\ell(x), &\text{if} \,\,\, 1 \le x \le \ell-1 \,\,\, \text{and} \,\,\, \si_\ell(x) \ne \ell, \\
               \xi(\ell-1), &\text{if} \,\,\, 1 \le x \le \ell-1 \,\,\, \text{and} \,\,\, \si_\ell(x) = \ell, \\
               \xi(x), &\text{if} \,\,\, \ell \le x \le k, \\
                   \end{cases} 
$$
$$
(\si\eta)_k(x) =   \begin{cases}
               \si_\ell(x), &\text{if} \,\,\, 1 \le x \le \ell-1, \\
               \eta(x), &\text{if} \,\,\, \ell \le x \le k, \qquad\qquad\qquad\qquad\quad \\
                   \end{cases} 
$$
and
$$
(\rho\eta)_k(x) =   \begin{cases}
               \rho_\ell(x), &\text{if} \,\,\, 1 \le x \le \ell-1, \\
               \eta(x), &\text{if} \,\,\, \ell \le x \le k \,\,\, \text{and} \,\,\, \eta(x) \ne \ell-1, \\
               \rho_\ell(\ell), &\text{if} \,\,\, \ell \le x \le k \,\,\, \text{and} \,\,\, \eta(x) = \ell-1, \quad \\
       \end{cases} 
$$
Then $(\si\xi)_k$ is a left extension of $\xi$ and $(\rho\eta)_k$ is a right extension of $\rho_\ell$ and the following hold:
\begin{itemize}
\item[\rm{(a)}]
If $\si_\ell$ and $\rho_\ell$ are right similar and $\xi$ and $\eta$ are two-sided similar, then $(\si\xi)_k$ and  $(\si\eta)_k$ are right and two-sided similar and $(\si\eta)_k$ and $(\rho\eta)_k$ are right similar.

\item[\rm{(b)}]
If $\si_\ell$ and $\rho_\ell$ are two-sided similar and $\xi$ and $\eta$ are two-sided similar, then $(\si\xi)_k$ and  $(\si\eta)_k$ are right and two-sided similar and $(\si\eta)_k$ and $(\rho\eta)_k$ are left and two-sided similar.
\end{itemize}}

\begin{figure}[htb]
\centerline{\epsfig{file=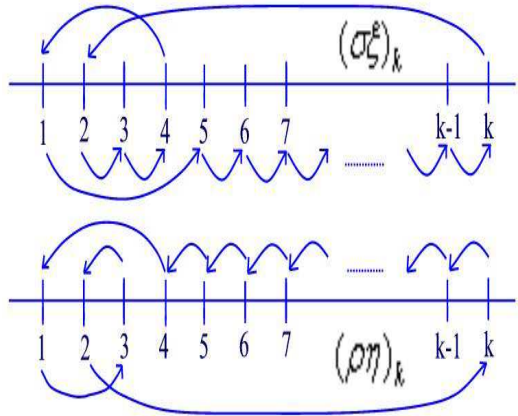,width=7cm,height=6cm}}
\caption{The cyclic permutations $(\si\xi)_k$ and $(\rho\eta)_k$.}
\end{figure}

\noindent
{\bf Remark.}
Later on, we will see that, for any $m \ge 6$, the cyclic permutations $\xi : 4 \ra m \ra m-1 \ra m-2 \ra \cdots \ra 6 \ra 5 \ra 4$ and $\eta : 4 \ra 5 \ra 6 \ra \cdots \ra m-1 \ra m \ra 4$ are two-sided similar (Theorem 8) and the cyclic permutations $\si_5 : 1 \ra 3 \ra 2 \ra 5 \ra 4 \ra 1$ and $\rho_5 : 1 \ra 5 \ra 2 \ra 3 \ra 4 \ra 1$ are right similar (Theorem 12).  Therefore, by Theorem 4(a), for any $k \ge 6$, the cyclic permutations $(\si\xi)_k : 1 \ra 3 \ra 2 \ra k \ra k-1 \ra k-2 \ra \cdots \ra 7 \ra 6 \ra 5 \ra 4 \ra 1$ and $(\rho\eta)_k : 1 \ra 5 \ra 6 \ra 7 \ra \cdots \ra k-1 \ra k \ra 2 \ra 3 \ra 4 \ra 1$ (see Figure 5) are right similar.  However, it can be easily checked that they are neither left nor two-sided similar (cf. Theorem 12). 
 
\noindent
{\it Proof.}
The desired results follow easily from the observations that the pair $((\si\xi)_k, (\si\eta)_k))$ is a synchronized left extension of $\xi$ and $\eta$ and the pair $((\si\eta)_k, (\rho\eta)_k)$ is a synchronized right extension of $\si_\ell$ and $\rho_\ell$.
\hfill\sq

Although two permutations with similar corresponding Petrie matrices need not be right, left, or two-sided, similar or weakly similar, we can still obtain more examples by "expanding" them to "longer" permutations.  

\noindent
{\bf Theorem 5.}
{\it Let $m \ge 0, n \ge 0$ and $k \ge 3$ be integers with $m+n \ge 1$.  Let $\sigma_k$ and $\rho_k$ be permutations on $P_k$ and let $\sigma_{m+k+n}$ and $\rho_{m+k+n}$ be permutations on $P_{m+k+n}$ such that
\begin{itemize}
\item[\rm {(1)}]
If $m = 0$, then $\sigma_{k+n}(i) = \sigma_k(i)$ and $\rho_{k+n}(i) = \rho_k(i)$ for all integers $1 \le i \le k$ and $\sigma_{k+n}(j) = \rho_{k+n}(j) \,\,\, (\ge k+1)$ for all integers $k+1 \le j \le k+n$,

\item[\rm {(2)}]
If $n = 0$, then $\sigma_{m+k}(j) = \rho_{m+k}(j) \,\,\, (\le m)$ for all integers $1 \le j \le m$ and $\sigma_{m+k}(j) = m+\sigma_k(j-m)$ and $\rho_{m+k}(j) = m+\rho_k(j-m)$ for all integers $m+1 \le j \le m+k$,

\item[\rm {(3)}]
If $m \ge 1$ and $n \ge 1$, then $\sigma_{m+k+n}(i) = \rho_{m+k+n}(i) \le m$ for all integers $1 \le i \le m$, and $\sigma_{m+k+n}(j) = m + \sigma_k(j-m)$ and $\rho_{m+k+n}(j) = m + \rho_k(j-m)$ for all integers $m+1 \le j \le m+k$ and $\sigma_{m+k+n}(x) = \rho_{m+k+n}(x) \ge m+k+1$ for all integers $m+k+1 \le x \le m+k+n$.
\end{itemize}

Assume that the Petrie matrices of $\sigma_k$ and $\rho_k$ are similar.  Then the following hold:
\begin{itemize}
\item[{\rm(1)}]
If both $\sigma_k$ and $\rho_k$ fix at least one same endpoint of $[1, k]$, then the following hold:
\begin{itemize}
\item[{\rm (a)}]
If $\sigma_k(k) = k = \rho_k(k)$, then $\sigma_k$ and $\rho_k$ are right similar.

\item[{\rm (b)}]
If $\sigma_k(1) = 1 = \rho_k(1)$, then $\sigma_k$ and $\rho_k$ are left similar.

\item[{\rm (c)}]
If $\sigma_k(k) = k = \rho_k(k)$ and $\sigma_k(1) = 1 = \rho_k(1)$, then $\sigma_k$ and $\rho_k$ are two-sided similar.
\end{itemize}

\item[{\rm (2)}]
If both $\sigma_k$ and $\rho_k$ fix no same endpoint of $[1, k]$, then the Petrie matrices of $\sigma_{m+k+n}$ and $\rho_{m+k+n}$ have the same characteristic polynomial and the following also hold:
\begin{itemize}
\item[{\rm (a)}]
If $m = 0$, then $\sigma_{k+n}$ and $\rho_{k+n}$ are right weakly similar.  Furthermore, if 1 is not an eigenvalue of the Petrie matrix of $\rho_k$, then the Petrie matrices of $\sigma_{k+n}$ and $\rho_{k+n}$ are similar and, $\sigma_{k+n}$ and $\rho_{k+n}$ are right similar.

\item[{\rm (b)}]
If $n = 0$, then $\sigma_{m+k}$ and $\rho_{m+k}$ are left weakly similar.  Furthermore, if 1 is not an eigenvalue of the Petrie matrix of $\rho_k$, then the Petrie matrices of $\sigma_{m+k}$ and $\rho_{m+k}$ are similar and, $\sigma_{m+k}$ and $\rho_{m+k}$ are left similar.

\item[{\rm (c)}]
If $m \ge 1$ and $n \ge 1$, then $\sigma_{m+k+n}$ and $\rho_{m+k+n}$ are right, left, and two-sided weakly similar.  Furthermore, if 1 is not an eigenvalue of the Petrie matrix of $\rho_k$, then the Petrie matrices of $\sigma_{m+k+n}$ and $\rho_{m+k+n}$ are similar and, $\sigma_{m+k+n}$ and $\rho_{m+k+n}$ are right, left, and two-sided similar.
\end{itemize}
\end{itemize}}

\noindent
{\bf Remark.}
In the above theorem, if 1 is an eigenvalue, then the result need not hold.  For example, let $\sigma = (13)$ and $\rho = (13)(24)$.  When considered as permutations on $P_4$, the Petrie matrix of $\sigma$ with respect to the basis $\{ J_1, J_2, J_3 \}$ is the same as that of $\rho$ with respect to the basis $\{ J_1, J_3, J_1+J_2+J_3 \}$.  So, the $3 \times 3$ Petrie matrices of $\sigma$ and $\rho$ are similar and since $\vp_{\sigma}(J_1+J_2) = J_1+J_2$, 1 is an eigenvalue.  However, the Petrie matrices of the permutations $(13)(56)$ and $(13)(24)(56)$ are not similar because they have distinct minimal polynomials.  Therefore, the permutations $(13)$ and $(13)(24)$, when considered as permutations on $P_5$, are not right similar although they are right weakly similar.  Their respective Petrie matrices are not similar either because they satisfy distinct minimal polynomials.  

\noindent
{\it Proof.}
Suppose $m = 0$.  Since the Petrie matrices of $\sigma_k$ and $\rho_k$ are similar, there is a basis $L = \{ L_1, L_2, \cdots, L_{k-1} \}$, where the $L_i$'s are linear combinations of $J_1, J_2, \cdots, J_{k-1}$ for $W_{{{\mathbb{R}}}^{k-1}}$ such that if $g$ is the linear transformation from $W_{{{\mathbb{R}}}^{k-1}}$ onto itself defined by putting $g(J_i) = L_i$ for all $1 \le i \le k-1$ then we have $(\vp_{\rho_k} \circ g)(x) = (g \circ \vp_{\sigma_k})(x)$ on $W_{{{\mathbb{R}}}^{k-1}}$.  Let $I_{\ell}$ denote the $\ell \times \ell$ identity matrix.  Then by expanding the determinants $\det(M_{\sigma_{k+n, k+n-1}} - \lambda I_{k+n-1})$ and $\det(M_{\rho_{k+n, k+n-1}} - \lambda I_{k+n-1})$ along the $(k+1)^{\rm {st}}$ column using Laplace's formula, we easily obtain that these two determinants are equal.  That is, the Petrie matrices of $\sigma_{k+n}$ and $\rho_{k+n}$ have the same characteristic polynomial.  Since $n \ge 1$ is arbitrary and $\sigma_{k+n} (= \rho_{k+n})$ is arbitrary on the set $\{ k+1, k+2, k+3, \cdots, k+n-1 \}$, this implies that $\sigma_{k+n}$ and $\rho_{k+n}$ are right weakly similar.    

On the other hand, if 1 is not an eigenvalue of the Petrie matrix of $\rho_k$, let ${\rm Id}_{k-1}$ be the identity map on $W_{{{\mathbb{R}}}^{k-1}}$.  Then the linear transformation $\vp_{\rho_k} - {\rm {Id}}_{k-1}$ is invertible.  Hence, there exists real numbers $a_1, a_2, \cdots, a_{k-1}$ such that $(\vp_{\rho_k} - {\rm {Id}}_{k-1})(\sum\limits_{i=1}^{k-1} a_iL_i) = \sum\limits_{i=\sigma_k(k)}^{k-1} L_i - \sum\limits_{i=\rho_k(k)}^{k-1} J_i$.  Let $S = \{ S_1, S_2, S_3, \cdots, S_{k+n-1} \} = \{ L_1, L_2, L_3, \cdots, L_{k-1}, \sum\limits_{i=1}^{k-1} a_iL_i + J_k, J_{k+1}, J_{k+2}, \cdots, J_{k+n-1} \}$ (in that order).  Then it is clear that $S$ is a basis for $W_{{{\mathbb{R}}}^{k+n-1}}$.  Let $h$ be the linear transformation from $W_{{{\mathbb{R}}}^{k+n-1}}$ onto itself defined by putting $h(J_i) = S_i$ for all $1 \le i \le k+n-1$.  Then we have $(\vp_{\rho_{k+n}} \circ h)(x) = (h \circ \vp_{\sigma_{k+n}})(x)$ on $W_{{{\mathbb{R}}}^{k+n-1}}$.  Therefore, the Petrie matrices of $\sigma_{k+n}$ and $\rho_{k+n}$ are similar.  Consequently, since any synchronized right extension of $\sigma_{k+n}$ and $\rho_{k+n}$ are just another different $\sigma_{k+n}$ and $\rho_{k+n}$, and so their respective Petrie matrices are similar.  This implies that $\sigma_{k+n}$ and $\rho_{k+n}$ are right similar.  This confirms (2a).  (2b) and (2c) can be proved similarly.  
\hfill\sq

\section{Examples of permutations which are right, left or two-sided, similar or weakly similar}
The following result is taken from {\bf{\cite{bl, go}}}.

\noindent
{\bf Theorem 6.}
{\it Let $n \ge 3$ be a fixed integer.  Let ${\mathbb{Z}}_2 = \{ 0, 1 \}$ denote the field of two elements.  Let $\si$ and $\rho$ be any two maps from $P_n$ into itself such that neither $\si$ nor $\rho$ nor $\rho \circ \si$ assumes the same value at any two consecutive integers $1 \le i < i+1 \le n$.  Then we have $\varphi_{\rho \circ \si} = \vpr \circ \vps$ on $W_{{{\mathbb{Z}_2}}^{n-1}}$.  In particular, if $\si$ is a permutation on $P_n$, then $\vps$ is an isomorphism on $W_{{{\mathbb{R}}}^{n-1}}$ and the Petrie matrix of $\si$ has determinant $\pm 1$.}

{\it Proof.}
For any two distinct real numbers $a$ and $b$, let $[a : b]$ denote the interval $[a, b]$ if $a < b$ or the interval $[b, a]$ if $a > b$.  Since the vector space $W_{{{\mathbb{Z}_2}}^{n-1}}$ is over the field ${\mathbb{Z}}_2$ of two elements, for any map $\gamma$ on $P_n$ with $\gamma(i) \ne \gamma(i+1)$ for all $1 \le i \le n-1$ and any two integers $1 \le k < \ell \le n$, we easily obtain that $\varphi_\gamma([k, \ell]) = [\gamma(k) : \gamma(\ell)]$.  Therefore, $\vpr(\vps(J_i)) = \vpr([\si(i) : \si(i+1)]) = [\rho(\si(i)) : \rho(\si(i+1))] = (\rho \circ \si)(J_i)$.  This shows that $\varphi_{\rho \circ \si} = \vpr \circ \vps$ on $W_{{{\mathbb{Z}_2}}^{n-1}}$.  
\hfill\sq

\begin{figure}[htb]
\centerline{\epsfig{file=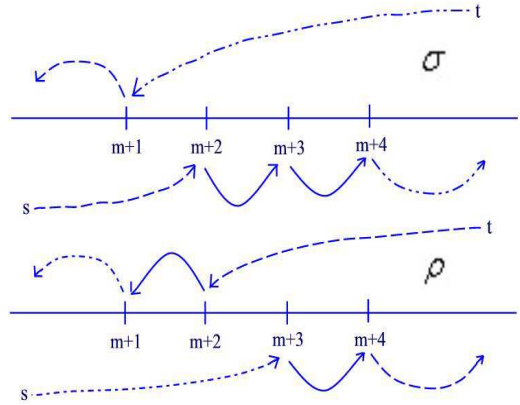,width=7cm,height=5.5cm}}
\caption{The permutations $\si$ and $\rho$ in Theorem 7.}
\end{figure}

We now come to the main result which initiates all of these (cf. {\bf{[8}}, Theorem 2{\bf ]}).

\noindent
{\bf Theorem 7.}
{\it Let $m \ge 1$ and $n \ge 0$ be integers.  Let $\si$ be a permutation on $P_{m+n+4}$ such that  
\begin{itemize}
\item[\rm{(1)}]
$\si(m+2) = m+3$,

\item[\rm{(2)}]
$\si(m+3) = m+4$,

\item[\rm{(3)}]
$\si(s) = m+2$ \, for some integer $s$ such that $1 \le s \le m$,

\item[\rm{(4)}]
$\si(i) \le m$ \, for all integers $i \ne s$ and $1 \le i \le m+1$,

\item[\rm{(5)}]
$\si(j) \ge m+5$ \, for all integers $j \ne t$ and $m+4 \le j \le m+4+n$ if $n \ge 1$,

\item[\rm{(6)}]
$\si(t) = m+1$ \, for some integer $t$ such that $m+5 \le t \le m+4+n$ if $n \ge 1$,

\item[\rm{(7)}]
$\si(m+4) = m+1$ \, if $n = 0$.
\end{itemize}

\noindent
We also let $\rho$ be a permutation on $P_{m+n+4}$ such that
\begin{itemize}
\item[\rm{(1)}]
$\rho(m+2) = m+1$, 

\item[\rm{(2)}]
$\rho(m+3) = m+4$,

\item[\rm{(3)}]
$\rho(s) = m+3$,

\item[\rm{(4)}]
$\rho(t) = m+2$ \, if $n \ge 1$,

\item[\rm{(5)}]
$\rho(m+4) = m+2$ \, if $n = 0$,

\item[\rm{(6)}]
$\rho(i) = \si(i)$ \, elsewhere.
\end{itemize}

\noindent
Then the Petrie matrices of the permutations $\si$ and $\rho$ are similar.}

\noindent
{\bf Remark.}
Note that our assumption implies that $\si$ and $\rho$ have exactly the same action (and hence are indistinguishable) on the set $\{ \, i : 1 \le i \le m+1$, $i \ne s \, \}$ and, if $n \ge 1$, also on the set $\{ \, j : m+4 \le j \le m+n+4$, $j \ne t \, \}$.  This motivates our definitions in section 2 of right, left and two-sided similarities of permutations on a finite set.  

\noindent
{\it Proof.}
Recall that, for every permutation $\ld$ on the set $\{ 1, 2, \cdots, m+n+1, m+n+2, m+n+3, m+n+4 \}$, the linearization $f_\ld$ of $\ld$ on the interval $[1, m+n+4]$ determines a linear transformation $\vpl$ on $W_{{{\mathbb{R}}}^{m+n+3}}$ such that, for real numbers $a_i$'s,
$$
\vpl \biggl(\sum\limits_{i=1}^{m+n+3} a_iJ_i \biggr) = \sum\limits_{i=1}^{m+n+3} 
a_i \vpl(J_i), \quad \text{where} \quad \vpl(J_i)= \sum\limits_{j=k_i}^{\ell_i} 
J_j \quad \text{if} \quad f_\ld(J_i) = \bigcup_{j=k_i}^{\ell_i} J_j.
$$

We now show that the set $S = \{ \, J_1, J_2, \cdots, J_m, J_{m+1}+J_{m+2}, J_{m+3}$, $\vpr(J_{m+3}), \vpr(J_{m+4})$, 
$\vpr(J_{m+5}), \cdots, \vpr(J_{m+n+3}) \, \}$ forms a basis for $W_{{{\mathbb{R}}}^{m+n+3}}$.  For $n \ge 0$, let $\om$ be the permutation on $P_{m+n+4}$ defined by putting $\om(i) = i$ for all integers $1 \le i \le m$, $\om(m+1) = m+3$, and $\om(j) = \rho(j)$ for all integers $m+2 \le j \le m+n+4$.  Then $\vpo(J_i) = J_i$ for all integers $1 \le i \le m-1$, $\vpo(J_m) = J_m + J_{m+1} + J_{m+2}$, $\vpo(J_{m+1}) = J_{m+1}+J_{m+2}$, $\vpo(J_{m+2}) = J_{m+1} + J_{m+2} + J_{m+3}$, and $\vpo(J_j) = \vpr(J_j)$ for all integers $m+3 \le j \le m+n+3$.  It is easy to see that the set $S$ is equal to the set $T = \{ \, \vpo(J_1), \vpo(J_2), \cdots, \vpo(J_{m-1}), \vpo(J_m) - \vpo(J_{m+1}), \vpo(J_{m+1})$, $\vpo(J_{m+2}) - \vpo(J_{m+1})$, $\vpo(J_{m+3})$, $\vpo(J_{m+4})$, $\vpo(J_{m+5}), \cdots$, $\vpo(J_{m+n+3}) \, \}$.  Therefore, by elementary operations, the determinant of the matrix $M(T\,|\,\{J_1, J_2, J_3, \cdots, J_{m+n+3}\})$ equals that of the matrix $M(\{ \vpo(J_i) : 1 \le i \le m+n+3 \}\, | \,\{J_1, J_2, \cdots, J_{m+n+3}\})$ which is, by Theorem 6, $\pm 1$.  Consequently, the determinant of the matrix $M(S\,|\,\{J_1, J_2, J_3, \cdots, J_{m+n+3}\})$ equals $\pm 1$.  This shows that $S$ is a basis for $W_{{{\mathbb{R}}}^{m+n+3}}$.

Let $h$ be the linear isomorphism from $W_{{{\mathbb{R}}}^{m+n+3}}$ onto itself defined by
$$
\begin{cases}
h(J_i) = J_i, \,\, 1 \le i \le m, \\
h(J_{m+1}+J_{m+2}) = J_{m+1}, \\
h(J_{m+3}) = J_{m+2}, \\
h(\vpr(J_i)) = J_i, \,\, m+3 \le i \le m+n+3. \\
\end{cases}
$$

We now show that $(h \circ \vpr)(L) = (\vps \circ h)(L)$ for every $L \in S$.  We have six cases to consider:

Case 1.  For some integers $1 \le i \le m$ and $1 \le k_i \le m$, $\vpr(J_i) = [k_i, m+3]$ (since $\rho(s) = m+3$, this means that $i = s-1$ or $s$).  In this case, we have $h(\vpr(J_i)) = h([k_i, m+3]) = h([k_i, m+1] + J_{m+1} + J_{m+2}) = h([k_i, m+1]) + h(J_{m+1} + J_{m+2}) = [k_i, m+1] + J_{m+1} = [k_i, m+2] = \vps(J_i) = \vps(h(J_i))$.  Thus, $(h \circ \vpr)(J_i) = (\vps \circ h)(J_i)$.

Case 2.  $1 \le i \le m$ and $m+3 \notin \vpr(J_i)$.  In this case, we have $\vpr(J_i) \sub [1, m+1]$.  So, $h(\vpr(J_i)) = \vpr(J_i) = \vps(J_i) = \vps(h(J_i))$.  Thus, $(h \circ \vpr)(J_i) = (\vps \circ h)(J_i)$.

Case 3.  $L = J_{m+1}+J_{m+2}$.  In this case, we have $h(\vpr(J_{m+1}+J_{m+2}) = h(\vpr(J_{m+1}) + \vpr(J_{m+2})) =  h([\rho(m+1), m+1] + J_{m+1} + J_{m+2} + J_{m+3}) = h([\rho(m+1), m+1]) + h(J_{m+1}+J_{m+2}) + h(J_{m+3}) = [\rho(m+1), m+1] + J_{m+1} + J_{m+2} = [\rho(m+1), m+3] = [\si(m+1), m+3] = \vps(J_{m+1}) = \vps(h(J_{m+1}+J_{m+2}))$.  Thus, $(h \circ \vpr)(J_{m+1}+J_{m+2}) = (\vps \circ h)(J_{m+1}+J_{m+2})$.

Case 4.  $L = J_{m+3}$.  In this case, we have $(h \circ \vpr)(J_{m+3}) = h(\vpr(J_{m+3})) = J_{m+3}$ (by definition of $h$)  $= \vps(J_{m+2}) = \vps(h(J_{m+3})) = (\vps \circ h)(J_{m+3})$.

Case 5.  $L = \vpr(J_i) = [m+2, \ell_i]$ for some integers $m+4 \le i \le m+n+3$ and $m+4 \le \ell_i \le m+n+4$ (since $\rho(t) = m+2$, this means that $i = t$ or $t-1$).  In this case, we have $h(\vpr(L)) = h(\vpr(\vpr(J_i))) = h(\vpr([m+2, \ell_i])) = h(\vpr(J_{m+2} + J_{m+3} + [m+4, \ell_i])) = h(J_{m+1}+J_{m+2}+J_{m+3} + \vpr(J_{m+3}) + \vpr([m+4, \ell_i])) = h(J_{m+1}+J_{m+2}) + h(J_{m+3}) + h(\vpr(J_{m+3})) + h(\vpr([m+4, \ell_i])) = J_{m+1} + J_{m+2} + J_{m+3} + [m+4, \ell_i] = [m+1, \ell_i] = \vps(J_i) = \vps(h(\vpr(J_i))) = \vps(h(L))$.  Thus, $(h \circ \vpr)(L) = (\vps \circ h)(L)$.

Case 6.  $L = \vpr(J_i)$ for some integer $m+4 \le i \le m+n+3$ and $m+2 \notin L$.  In this case, $\vpr(J_i) \sub [m+4, 
m+n+4]$.  So, $\vpr(L) = \vpr(\vpr(J_i)) \sub \vpr([m+4, m+n+4])$.  Thus, $h(\vpr(L)) = h(\vpr(\vpr(J_i))) = \vpr(J_i)
= \vps(J_i) = \vps(h(\vpr(J_i))) = \vps(h(L))$.  Hence, $(h \circ \vpr)(L) = (\vps \circ h)(L)$.

Therefore, we have shown that $(h \circ \vpr)(L) = (\vps \circ h)(L)$ for every $L \in S$.  Thus, $h \circ \vpr =
\vps \circ h$ on $W_{{{\mathbb{R}}}^{m+n+3}}$.  Since $h$ is an isomorphism from $W_{{{\mathbb{R}}}^{m+n+3}}$ onto itself, this implies that $\vpr$ is conjugate to $\vps$ through $h$.  Hence the Petrie matrix of $\si$ is similar {\bf{\cite{he}}} to that of $\rho$.  
\hfill\sq
\bigskip

The following result is an easy consequence of Theorem 7.

\noindent
{\bf Theorem 8.}
{\it Let $k \ge 4$, and $3 \le n \le k-1$ be integers.  Let $\si_{n,k}$ be the cyclic permutation from $P_k$ onto itself such that
\begin{itemize}
\item[\rm{(1)}]
$\si_{n,k}(1) = n$,

\item[\rm{(2)}]
$\si_{n,k}(i) = i+1$ \, for all \, $n \le i \le k-1$,

\item[\rm{(3)}]
$\si_{n,k}(k) = n-1$

\item[\rm{(4)}]
$\si_{n,k}(j) = j-1$ \, for all \, $2 \le j \le n-1$.
\end{itemize}

\noindent
Furthermore, let $\si_{2,k}$ denote the cyclic permutation on $P_k$ such that $\si_{2,k}(i) = i+1$ for all $1 \le i \le k-1$ and $\si_{2,k}(k) = 1$, and let $\si_{k,k}$ denote the cyclic permutation on $P_k$ such that $\si_{k,k}(1) = k$ and $\si_{k,k}(j) = j-1$ for all $2 \le j \le k$.  Then the following hold:
\begin{itemize}
\item[\rm{(1)}]

If $k \ge 5$, then all $\si_{n,k}$, $3 \le n \le k-1$, are right, left, and two-sided similar to one another.  Moreover, the Petrie matrices of $\si_{n,k}$, $3 \le n \le k-1$, are similar to one another.

\item[\rm{(2)}]
If $k \ge 4$ and $n$ is any integer with $3 \le n \le k-1$, then $\si_{2,k}$ and $\si_{n,k}$ are left and two-sided, but not right, similar.  Moreover, the Petrie matrices of $\si_{2,k}$ and $\si_{n,k}$ are not similar because they have distinct traces.

\item[\rm{(3)}]
If $k \ge 4$ and $n$ is any integer with $3 \le n \le k-1$, then $\si_{n,k}$ and $\si_{k,k}$ are right and two-sided, but not left, similar.  Moreover, the Petrie matrices of $\si_{n,k}$ and $\si_{k,k}$ are not similar because they have distinct traces.

\item[\rm{(4)}]
If $k \ge 4$, then $\si_{2,k}$ and $\si_{k,k}$ are two-sided, but neither right nor left, similar.  Moreover, the Petrie matrices of $\si_{2,k}$ and $\si_{k,k}$ are similar.
\end{itemize}}

\noindent
{\bf Remark.}
It follows from (4) of the above result that, for any two (cyclic) permutations on $P_k$, even their respective Petrie matrices are similar, there is no guarantee that they are one-sided (right or left) similar.  They need not be two-sided similar either, see examples in Theorem 12 below.

\noindent
{\it Proof.}
For each integer $3 \le i \le n-2$, the Petrie matrices of $\si_{i,k}$ and $\si_{i+1,k}$ are similar by Theorem 7.  So, the Petrie matrices of $\si_{n,k}$, $3 \le n \le k-1$, are similar to one another.  Furthermore, it follows from Theorem 6 that, for each integer $2 \le i \le k-2$, the cyclic permutations $\si_{i,k}$ and $\si_{i+1,k}$ are left and two-sided similar.  Therefore, all $\si_{n,k}$, $2 \le n \le k-1$, are left and two-sided similar to one another.  Similarly, by Theorem 7, all $\si_{n,k}$, $3 \le n \le k$, are right and two-sided similar to one another.  In particular, all $\si_{n,k}$, $3 \le n \le k-1$, are right, left and two-sided similar to one another.  

On the other hand, by direct computations, for any integers $k > n \ge 3$, the trace of the Petrie matrix of the cyclic permutation $1 \ra 2 \ra 3 \ra \cdots \ra k-1 \ra k \ra k+1 \ra 1$ is $1$ while the trace of the Petrie matrix of the cyclic permutation $1 \ra n \ra n+1 \ra \cdots \ra k-1 \ra k \ra k+1 \ra n-1 \ra n-2 \ra \cdots \ra 3 \ra 2 \ra 1$ is $3$.  So, for any $3 \le n \le k-1$, the Petrie matrices of the cyclic permutations $1 \ra 2 \ra 3 \ra \cdots \ra k-1 \ra k \ra k+1 \ra 1$ and $1 \ra n \ra n+1 \ra \cdots \ra k-1 \ra k \ra k+1 \ra n-1 \ra n-2 \ra \cdots \ra 3 \ra 2 \ra 1$ are not similar.  Therefore, for any $3 \le n \le k-1$, the cyclic permutations $\si_{2,k}$ and $\si_{n,k}$ are not right similar and the Petrie matrices of $\si_{2,k}$ and $\si_{n,k}$ are not similar.  Similarly, for any $3 \le n \le k-1$, the cyclic permutations $\si_{k,k}$ and $\si_{n,k}$ are not left similar and the Petrie matrices of $\si_{k,k}$ and $\si_{n,k}$ are not similar.  
\hfill\sq
\bigskip

We shall need the following result.

\noindent
{\bf Lemma 9.} 
Let $k \ge 3$ and $n \ge 3$ be fixed integers and let $\mu$ be any permutation on $P_{k+n}$ such that $\mu(i) = i-1$ for all integers $2 \le i \le k-1$ and $\mu(1) = k < \mu(k)$.  Then the set $S = \biggl\{ \, \sum\limits_{i=1}^{k-1} J_i, [k, \mu(k)]$, $\vpm([k, \mu(k)])$, $\vpm^2([k, \mu(k)])$, $\vpm^3([k, \mu(k)])$, $\cdots, \vpm^{k-3}([k, \mu(k)])$, $\vpm^{k-2}(J_k)$, $\vpm^{k-2}(J_{k+1})$, $\vpm^{k-2}(J_{k+2})$, $\cdots, \vpm^{k-2}(J_{k+n-1}) \biggr\}$ equals the set $\{ \, \vpm^{k-2}(J_i) : 1 \le i \le k+n-1 \, \}$ and so is a basis for $W_{{{\mathbb{R}}}^{k+n-1}}$.  Furthermore, the determinant of $M(S|B_{n+k-1})$ is $\pm 1$, where $B_{n+k-1}$ is the standard basis for $W_{{{\mathbb{R}}}^{n+k-1}}$.

\noindent
{\it Proof.}
By direct computations, we obtain that $\sum\limits_{\ell=1}^{k-1} J_\ell = \vpm^{k-2}(J_{k-2})$, and, since $[k, \mu(k)] = \vpm^{k-2}(J_{k-3}) - 2\vpm^{k-2}(J_{k-2})$, we have, for every integer $0 \le i \le k-4$, $$\vpm^i([k, \mu(k)]) = \vpm^{k-2}(\vpm^i(J_{k-3})) - 2\vpm^{k-2}(\vpm^i(J_{k-2})) = \vpm^{k-2}(J_{k-i-3}) - 2\vpm^{k-2}(J_{k-i-2}),$$ and $$\sum\limits_{i=2}^{k-1} \vpm^{k-2}(J_i) - \vpm^{k-2}(J_1) = \vpm^{k-3}((\sum\limits_{i=2}^{k-1} \vpm(J_i)) - \vpm(J_1)) = \vpm^{k-3}([k, \mu(k)]).$$  Therefore, $S = \{ \, \vpm^{k-2}(J_{k-2})$, $\vpm^{k-2}(J_{k-3}) - 2\vpm^{k-2}(J_{k-2})$, $\vpm^{k-2}(J_{k-4})$ $-$ $2\vpm^{k-2}(J_{k-3})$, $\vpm^{k-2}(J_{k-5})$ $-$ $2\vpm^{k-2}(J_{k-4})$, $\cdots$, $\vpm^{k-2}(J_1)$ $-$ $2\vpm^{k-2}(J_2)$, $\sum\limits_{i=2}^{k-1} \vpm^{k-2}(J_i)$ $-$ $\vpm^{k-2}(J_1)$, $\vpm^{k-2}(J_k)$, $\vpm^{k-2}(J_{k+1})$, $\vpm^{k-2}(J_{k+2})$, $\cdots$, $\vpm^{k-2}(J_{k+n-1}) \, \}$.  By applying elementary operations to the set $S$, we obtain the set $T = \{ \, \vpm^{k-2}(J_i) : 1 \le i \le k+n-1 \, \}$.  By Theorem 6, the determinant of $M(T|B_{n+k-1})$ equals $(\pm 1)^{k-2}$.  Therefore, the determinant of $M(S|B_{n+k-1})$ equals $\pm 1$.  Consequently, $S$ is also a basis for $W_{{{\mathbb{R}}}^{k+n-1}}$.
\hfill\sq

\begin{figure}[htb]
\centerline{\epsfig{file=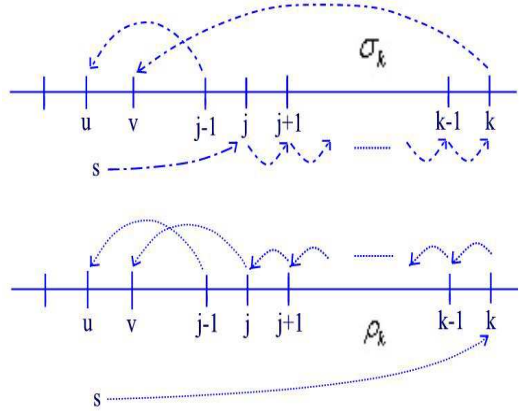,width=7cm,height=5.5cm}}
\caption{The permutations $\si_k$ and $\rho_k$ in Theorem 10.}
\end{figure}

The following result is similar to Theorem 4 in some sense.  

\noindent
{\bf Theorem 10.}
{\it Let $k > j \ge 3$ be integers and let $\si_k$ and $\rho_k$ be two permutations on $P_k$ such that 
\begin{itemize}
\item[(a)] There exists an integer $1 \le s \le j-2$ such that $\si_k(s) = j$ and $\rho_k(s) = k$;

\item[(b)] $\si_k(x) = \rho_k(x) \le j-1$ for all integers $1 \le x \le j-1$ and $x \ne s$;

\item[(c)] $\si_k(x) = x+1$ for all integers $j \le x \le k-1$;

\item[(d)] $\rho_k(x) = x-1$ for all integers $j+1 \le x \le k$;

\item[(e)] $\si_k(j-1) = \rho_k(j-1) < \si_k(k) = \rho_k(j) \le j-1$.  
\end{itemize}

Then $\si_k$ and $\rho_k$ are right and two-sided, but not left, similar.  Moreover, the Petrie matrices of $\si_k$ and $\rho_k$ are not similar because they have distinct traces.}

\noindent
{\bf Remarks.}
(1) Note that $\si_k$ and $\rho_k$ involve the following three permutations : The cyclic permutations $\al_{k-j+1} : 1 \ra 2 \ra 3 \ra \cdots \ra k-j \ra k-j+1 \ra 1$ and $\beta_{k-j+1} : 1 \ra k-j+1 \ra k-j \ra \cdots \ra 3 \ra 2 \ra 1$, and the permutation $\pi_j$, where $\pi_j(i) = \si_k(i)$ for all $1 \le i \le j-1$ and $\pi_k(j) = \si_k(k)$.  Here, $\si_k$ is a right extension of $\pi_j$ and $\rho_k$ is a left extension of $\beta_{k-j+1}$ and, by Theorem 8(4), $\al_{k-j+1}$ and $\beta_{k-j+1}$ are two-sided similar.  

(2) If $\si_k(k) = \rho_k(j) \le \si_k(j-1) = \rho_k(j-1) \le j-1$, then $\si_k$ and $\rho_k$ need not be right or two-sided similar.  For example, if $j = 4$ and $k \ge 5$ and $$\si_k : 1 \ra 3 \ra 2 \ra 4 \ra 5 \ra 6 \ra \cdots \ra k-1 \ra k \ra 1$$ and $$\rho_k : 1 \ra 3 \ra 2 \ra k \ra k-1 \ra k-2 \ra k-3 \ra \cdots \ra 6 \ra 5 \ra 4 \ra 1,$$ then the Petrie matrices of the synchronized right extensions $$1 \ra 3 \ra 2 \ra 4 \ra 5 \ra 6 \ra \cdots \ra k-1 \ra k \ra k+1 \ra 1$$ and $$1 \ra 3 \ra 2 \ra k \ra k+1 \ra k-1 \ra k-2 \ra k-3 \ra \cdots \ra 6 \ra 5 \ra 4 \ra 1$$ have distinct determinants $(-1)^{k+1}, (-1)^k$ respectively.  So, $\si_k$ and $\rho_k$ are not right similar.  Similarly, the Petrie matrices of the synchronized two-sided extensions : $$1 \ra 4 \ra 3 \ra 5 \ra 6 \ra 7 \ra \cdots \ra k-1 \ra k \ra k+1 \ra k+2 \ra 2 \ra 1$$ and $$1 \ra 4 \ra 3 \ra k+1 \ra k+2 \ra k \ra k-1 \ra k-2 \ra \cdots \ra 7 \ra 6 \ra 5 \ra 2 \ra 1$$ have distinct determinants $(-1)^k, (-1)^{k+1}$ respectively.  So, $\si_k$ and $\rho_k$ are not two-sided similar.  

\noindent
{\it Proof.}
For any positive integer $n$, let $(\si_{k+n}, \rho_{k+n})$ be a synchronized right extension of $\si_k$ and $\rho_k$.  Let $\vp_{\si_{k+n}}$ and $\rho_{\si_{k+n}}$ be the linear transformations on $W_{{{\mathbb{R}}}^{k+n-1}}$ determined by the linearizations of $\si_{k+n}$ and $\rho_{k+n}$ respectively as introduced in section 1.  Let $S = \{ \, J_1, J_2, \cdots, J_{j-2}, \sum\limits_{i=j-1}^{k-1}J_i$, $[k, \rho_{k+n}(k)]$, $\vp_{\rho_{k+n}}([k, \rho_{k+n}(k)])$, $\vp_{\rho_{k+n}}^2([k, \rho_{k+n}(k)])$, $\vp_{\rho_{k+n}}^3([k, \rho_{k+n}(k)], \cdots$, $\vp_{\rho_{k+n}}^{k-j-1}([k, \rho_{k+n}(k)]$, $\vp_{\rho_{k+n}}^{k-j}(J_{k})$, $\vp_{\rho_{k+n}}^{k-j}(J_{k+1})$, $\cdots$, $\vp_{\rho_{k+n}}^{k-j}(J_{k+n-1}) \, \}$.  Let $\mu$ be the permutation on the set $\{ j-1, j, j+1, j+2, \cdots, k+n-1, k+n \}$ such that $\mu(i) = i-1$ for all $j \le i \le k-1$, $\mu(j-1) = k$, and $\mu(i) = \rho_{k+n}(i)$ for all $k \le i \le k+n$.  By Lemma 9, the set $T = \{ \sum\limits_{i=j-1}^{k-1}J_i$, $[k, \rho_{k+n}(k)]$, $\vp_{\rho_{k+n}}([k, \rho_{k+n}(k)])$, $\vp_{\rho_{k+n}}^2([k, \rho_{k+n}(k)])$, $\cdots, \vp_{\rho_{k+n}}^{k-j-1}([k, \rho_{k+n}(k)]$, $\vp_{\rho_{k+n}}^{k-j}(J_{k})$, $\vp_{\rho_{k+n}}^{k-j}(J_{k+1})$, $\cdots, \vp_{\rho_{k+n}}^{k-j}(J_{k+n-1}) \, \}$ equals the set $\{ \vp_{\mu}^{k-j}(J_i) : j-1 \le i \le k+n-1 \}$.  So, by Theorem 6 and {\bf{\cite{he}}}, $T$ is a basis and so, $S$ is a basis for $W_{{{\mathbb{R}}}^{k+n-1}}$.  

Let $h$ be the linear isomorphism from $W_{{{\mathbb{R}}}^{k+n-1}}$ onto itself defined by
$$
\begin{cases}
h(J_i) = J_i, \, \text{for} \, 1 \le i \le j-1, \\
h(J_{j-1}) = \sum\limits_{i=j-1}^{k-1}J_i, \\
h(J_{j+i}) = \vp_{\rho_{k+n}}^i([k, \rho_{k+n}(k)]), \,\, 0 \le i \le k-j-1, \\
h(J_i) = \vp_{\rho_{k+n}}^{k-j}(J_j), \,\, k \le i \le k+n-1. 
\end{cases}
$$
Then it is easy to see that $(h \circ \vp_{\si_{k+n}})(J_i) = (\vp_{\rho_{k+n}} \circ h)(J_i)$ for all $1 \le i \le k+n-1$.  Therefore, the Petrie matrix of $\si_{k+n}$ is similar to that of  $\rho_{k+n}$.  This shows that $\si_k$ and $\rho_k$ are right similar and, consequently, the two-sided similarity of $\si_k$ and $\rho_k$ also follows.

Let $\si_k$ and $\rho_k$ denote the cyclic permutations $1 \ra 4 \ra 5 \ra 6 \ra \cdots \ra k-1 \ra k \ra 2 \ra 3 \ra 1$ and $1 \ra k \ra k-1 \ra k-2 \ra \cdots \ra 5 \ra 4 \ra 2 \ra 3 \ra 1$ \, respectively.  The Petrie matrices of the cyclic permutations $1 \ra 5 \ra 6 \ra 7 \ra \cdots \ra k-1 \ra k \ra k+1 \ra 3 \ra 4 \ra 2 \ra 1$ and $1 \ra k+1 \ra k \ra k-1 
\ra \cdots \ra 6 \ra 5 \ra 3 \ra 4 \ra 2 \ra 1$ have distinct traces.  Therefore, $\si_k$ and $\rho_k$ are not left similar.
\hfill\sq

The above result presents a method to construct permutations on $P_k$ with $k \ge 5$ which are right and two-sided similar.  This, combined with Theorem 4, will produce even more examples which are right and two-sided similar.  In the following, we introduce some concrete examples by applying Theorem 4 to Theorems 8(4) and 10.

\begin{figure}[htb]
\centerline{\epsfig{file=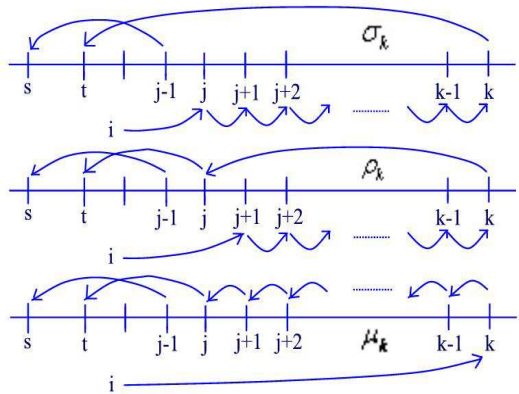,width=7cm,height=5.5cm}}
\caption{The permutations $\si_k$, $\rho_k$ and $\mu_k$ in Corollary 11.}
\end{figure}

\noindent
{\bf Corollary 11.}
{\it Let $j \ge 3$ be an integer and let $\pi_j$ be any fixed permutation on $P_j$ such that $\pi_j(j-1) < \pi_j(j) < j$.  For any integer $k \ge j+2$, let $\si_k$, $\rho_k$, $\mu_k$, and $\nu_k$ be the permutations on $P_k$ defined by 
$$
\si_k(x) =   \begin{cases}
               \pi_j(x), &\text{if} \,\,\, 1 \le x \le j-1, \qquad\qquad\quad\cr
               x+1, &\text{if} \,\,\, j \le x \le k-1, \\
               \pi_j(j), &\text{if} \,\,\, x = k, \\
       \end{cases} 
$$
$$
\rho_k(x) =  \begin{cases}
               \pi_j(x), &\text{if} \,\,\, 1 \le x \le j \,\, \text{and} \,\,\, \pi_j(x) \ne j, \cr
               j+1, &\text{if} \,\,\, 1 \le x \le j \,\, \text{and} \,\,\, \pi_j(x) = j, \\
               x+1, &\text{if} \,\,\, j < x \le k-1, \\
               j, &\text{if} \,\,\, x = k, \\
       \end{cases}
$$
$$
\mu_k(x) =  \begin{cases}
               \pi_j(x), &\text{if} \,\,\, 1 \le x \le j, \,\, \text{and} \,\,\, \pi_j(x) \ne j, \cr
               k, &\text{if} \,\,\, 1 \le x \le j, \,\, \text{and} \,\,\, \pi_j(x) = j, \\
               x-1, &\text{if} \,\,\, j < x \le k, \\
       \end{cases}
$$
\noindent
and
$$
\nu_k(x) =  \begin{cases}
               \pi_j(x), &\text{if} \,\,\, 1 \le x \le j-1, \qquad\qquad\quad \\
               k, &\text{if} \,\,\, x = j, \\
               x-1, &\text{if} \,\,\, j+2 \le x \le k, \\
               \pi(j), &\text{if} \,\,\, x = j+1. \\
       \end{cases}
$$
\noindent
Then $\si_k$, $\rho_k$ and $\mu_k$ are right and two-sided similar to one another, $\si_k$ and $\rho_k$ are also left similar while $\rho_k$ and $\mu_k$ need not be, and the Petrie matrices of $\si_k$ and $\rho_k$ are similar (because $\si_{j+1}$ and $\rho_{j+1}$ are right similar by Theorem 10) while those of $\rho_k$ and $\mu_k$ are not because they have distinct traces.  Furthermore, $\nu_k$ is neither right nor two-sided similar to any of $\si_k$, $\rho_k$, and $\mu_k$.}

\noindent
{\bf Remarks.}
(1) Note that $\si_k$ is a right extension of the permutation $\pi_j$, $\rho_k$ is a left extension of the cyclic permutation $j \ra j+1 \ra j+2 \ra \cdots \ra k-1 \ra k \ra j$, $\mu_k$ is a left extension of the cyclic permutation $j \ra k \ra k-1 \ra k-2 \ra \cdots \ra j+2 \ra j+1 \ra j$, and $\nu_k$ is another right extension of $\pi_j$ which is different from $\si_k$.  The above result says that, among these four various extensions, three are right and two-sided similar to one another while the remaining one is neither right nor two-sided similar to any of the three.  

\begin{figure}[htb]
\centerline{\epsfig{file=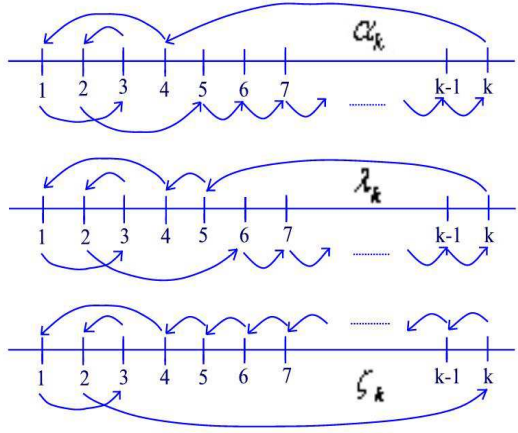,width=7cm,height=6cm}}
\caption{The permutations $\al_k$, $\ld_k$ and $\zeta_k$ in Remark(2) following Corollary 11.}
\end{figure}

(2) In the above corollary, if we take $j = 5$ and the permutation $\pi_5 : 1 \ra 3 \ra 2 \ra 5 \ra 4 \ra 1$, then the corollary says that $\al_k$, $\ld_k$ and $\zeta_k$ (see Figure 9), for any $k \ge 7$, are all right and two-sided similar to one another.  Also, by Theorem 10, $\al_6$ and $\ld_6$ are right and two-sided similar.  As for left similarity, it follows from Theorem 7 that $\al_k$ and $\ld_k$ are left similar for all $k \ge 7$ while $\al_6$ and $\ld_6$ are not left similar because of distinct traces for some easy synchronized left extensions.  Moreover, for each $k \ge 7$, $\zeta_k$ is neither left similar to $\al_k$ nor to $\ld_k$ because of distinct traces for some easy synchronized left extensions.  

\begin{figure}[htb]
\centerline{\epsfig{file=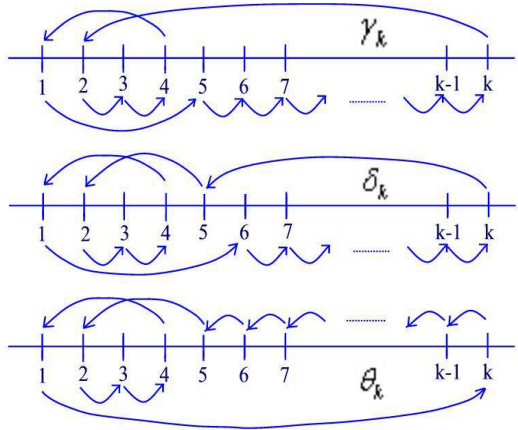,width=7cm,height=6cm}}
\caption{The permutations $\gamma_k$, $\delta_k$ and $\tha_k$ in Remark(3) following Corollary 11.}
\end{figure}

(3) In the above corollary, if we take $j = 5$ and the permutation $\pi_5 : 1 \ra 5 \ra 2 \ra 3 \ra 4 \ra 1$, then the corollary says that $\gamma_k$, $\delta_k$ and $\tha_k$(see Figure 10), for any $k \ge 7$, are all right and two-sided similar to one another.  On the other hand, by Theorem 10, $\gamma_6$ and $\delta_6$ are (right and) two-sided similar, and so, for any $k \ge 7$, $\gamma_k$ and $\delta_k$ are also {\it left} similar.  But, $\gamma_6$ and $\delta_6$ are not left similar because of the distinct traces of some easy synchronized left extensions.  Furthermore, for any $k \ge 7$, $\tha_k$  cannot be left similar to any of $\gamma_k$ and $\delta_k$ because of the distinct traces of some easy synchronized left extensions.    

\noindent
{\it Proof.}
By Theorem 10, $\si_{j+1}$ and $\rho_{j+1}$ are right and two-sided similar.  On the other hand, for $k \ge j+2$, $(\si_k, \rho_k)$ is a synchronized right extension of $\si_{j+1}$ and $\rho_{j+1}$.  This, combined with Theorem 2, implies that $\si_k$ and $\rho_k$ are right, left, and two-sided similar.  On the other hand, by Theorem 8(4), the cyclic permutations $\al_{k-j+1} : 1 \ra 2 \ra 3 \ra \cdots \ra k-j+1 \ra 1$ and $\be_{k-j+1} : 1 \ra k-j+1 \ra k-j \ra k-j-1 \ra \cdots \ra 3 \ra 2 \ra 1$ on $P_{k-j+1}$ are two-sided similar.  Since $(\rho_k, \mu_k)$ is a synchronized left extension of $\al_{k-j+1}$ and $\be_{k-j+1}$, we obtain that $\rho_k$ and $\mu_k$ are right and two-sided similar.  
\hfill\sq
\bigskip

In the following, we present more examples of {\it cyclic} permutations which are right, but neither left nor two-sided, similar.

\begin{figure}[htb]
\centerline{\epsfig{file=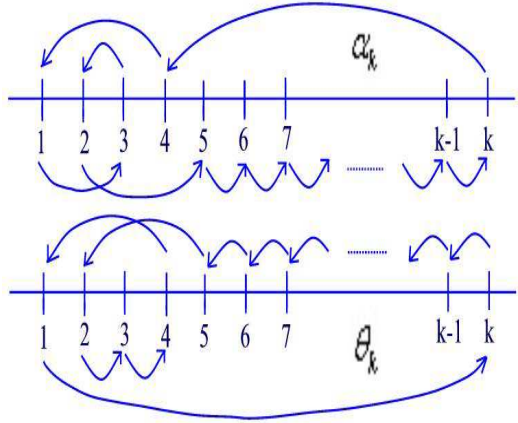,width=7cm,height=6cm}}
\caption{The cyclic permutations $\al_k$ and $\tha_k$ in Theorem 12.}
\end{figure}

\noindent
{\bf Theorem 12.}
Let $k \ge 5$ be an integer.  Then the cyclic permutations \\
\indent\indent\indent $\al_k : 1 \ra 3 \ra 2 \ra 5 \ra 6 \ra 7 \ra \cdots \ra k-2 \ra k-1 \ra k \ra 4 \ra 1$ and \\
\indent\indent\indent $\tha_k : 1 \ra k \ra k-1 \ra k-2 \ra \cdots \ra 7 \ra 6 \ra 5 \ra 2 \ra 3 \ra 4 \ra 1$ \\
on $P_k$ are right, but neither left nor two-sided, similar.  Moreover, their respective Petrie matrices are similar if $k = 5$ and not if $k > 5$. 

\noindent
{\bf Remark.}
It follows from the above result that the Petrie matrices of $\al_5$ and $\tha_5$ are similar.  However, they are not two-sided similar.

\noindent
{\it Proof.}
For any integer $n \ge 0$, let $(\al_{k+n}, \tha_{k+n})$ be a synchronized right extension of $\al_k$ and $\tha_k$ on $P_{k+n}$ and let $\vp_{\al_{k+n}}$ and $\vp_{\tha_{k+n}}$ be the linear transformations on $W_{{{\mathbb{R}}}^{k+n-1}}$  determined by the linearisations of $\al_{k+n}$ and $\tha_{k+n}$ respectively as introduced in section 1.  Let $S =  \{ \, J_1+J_2, J_1+J_2+J_3, J_3-J_2, \sum\limits_{i=2}^{k-1} J_i, [k, \tha_{k+n}(k)]$, $\vp_{\tha_{k+n}}([k, \tha_{k+n}(k)])$, $\vp_{\tha_{k+n}}^2([k, \tha_{k+n}(k)])$, $\cdots$, $\vp_{\tha_{k+n}}^{k-6}([k, \tha_{k+n}(k)])$, $\vp_{\tha_{k+n}}^{k-5}(J_k)$, $\vp_{\tha_{k+n}}^{k-5}(J_{k+1})$, $\vp_{\tha_{k+n}}^{k-5}(J_{k+2})$, $\cdots, \vp_{\tha_{k+n}}^{k-5}(J_{k+n-1}) \}$.  We want to show that $S$ is a basis for $W_{{{\mathbb{R}}}^{k+n-1}}$.  Let $V$ be the subspace spanned by $S$.  It is easy to see that $\{ J_1, J_2, J_3 \} \subset V$.  Furthermore, $$\vp_{\tha_{k+n}}(J_1) = \sum\limits_{i=2}^{k-1} J_i - J_2 \quad \in \quad < J_1, J_2, J_3, \sum\limits_{i=2}^{k-1} J_i, [k, \tha_{k+n}(k)] >,$$ $$\vp_{\tha_{k+n}}(J_2) = J_3 \quad \in \quad < J_1, J_2, J_3, \sum\limits_{i=2}^{k-1} J_i, [k, \tha_{k+n}(k)] >,$$ $$\vp_{\tha_{k+n}}(J_3) = J_1 + J_2 + J_3 \quad \in \quad < J_1, J_2, J_3, \sum\limits_{i=2}^{k-1} J_i, [k, \tha_{k+n}(k)] >,$$ $$\vp_{\tha_{k+n}}(J_4) = J_1 \quad \in \quad < J_1, J_2, J_3, \sum\limits_{i=2}^{k-1} J_i, [k, \tha_{k+n}(k)] >,$$ and $$\vp_{\tha_{k+n}}(\sum\limits_{i=2}^{k-1} J_i) = J_1 + J_3 + (J_1+J_2+J_3) + \sum\limits_{i=2}^{k-1} J_i + [k, \tha_{k+n}(k)] \,\, \in \,\, < J_1, J_2, J_3, \sum\limits_{i=2}^{k-1} J_i, [k, \tha_{k+n}(k)] >.$$  Therefore, for $1 \le i \le 4$, $$\vp_{\tha_{k+n}}^{k-5}(J_i) \,\, \in \,\,\, < J_1, J_2, J_3, \sum\limits_{i=2}^{k-1} J_i, [k, \tha_{k+n}(k)], \vp_{\tha_{k+n}}([k, \tha_{k+n}(k)]), \cdots, \vp_{\tha_{k+n}}^{k-6}([k, \tha_{k+n}(k)]) >.$$  On the other hand, for $2 \le i \le k-5$, $\vp_{\tha_{k+n}}^{k-5-i}(J_{k-i}) = J_5$.  So, $\vp_{\tha_{k+n}}^{k-3-i}(J_{k-i}) = \vp_{\tha_{k+n}}^2(J_5) =  \vp_{\tha_{k+n}}(J_2 + J_3 + J_4) \quad \in \quad < J_1, J_2, J_3, \sum\limits_{i=2}^{k-1} J_i, [k, \tha_{k+n}(k)] >$.  Consequently, $$\vp_{\tha_{k+n}}^{k-5}(J_{k-i}) \,\, \in \,\, < J_1, J_2, J_3, \sum\limits_{i=2}^{k-1} J_i, [k, \tha_{k+n}(k)], \vp_{\tha_{k+n}}([k, \tha_{k+n}(k)]), \cdots, \vp_{\tha_{k+n}}^{i-2}([k, \tha_{k+n}(k)]) >.$$  Finally, since $\vp_{\tha_{k+n}}(J_{k-1}) = J_{k-2}+ J_{k-1} + [k, \tha_{k+n}(k)]$, \, $\vp_{\tha_{k+n}}^{k-5}(J_{k-1}) = (J_2+J_3+J_4) + (J_5+J_6+ \cdots + J_{k-1} + [k, \tha_{k+n}(k)]+ \vp_{\tha_{k+n}}([k, \tha_{k+n}(k)])+ \cdots + \vp_{\tha_{k+n}}^{k-6}([k, \tha_{k+n}(k)]))$.  Therefore, for $5 \le i \le k-1$, we have $$\vp_{\tha_{k+n}}^{k-5}(J_{i}) \,\, \in \,\, < J_1, J_2, J_3, \sum\limits_{i=2}^{k-1} J_i, [k, \tha_{k+n}(k)], \vp_{\tha_{k+n}}([k, \tha_{k+n}(k)]), \cdots, \vp_{\tha_{k+n}}^{k-6}([k, \tha_{k+n}(k)])>.$$  This shows that the set $T = \{ \vp_{\tha_{k+n}}^{k-5}(J_i) : 1 \le i \le k+n-1 \}$ is contained in $V$.  Since, by Theorem 6, $T$ is a basis for $W_{{{\mathbb{R}}}^{k+n-1}}$, we obtain that $V = W_{{{\mathbb{R}}}^{k+n-1}}$.  That is, $S$ is a basis for $W_{{{\mathbb{R}}}^{k+n-1}}$.  

Let $h$ be the linear isomorphism from $W_{{{\mathbb{R}}}^{k+n-1}}$ onto itself defined by
$$
\begin{cases}
h(J_1) = J_1 + J_2, \\
h(J_2) = J_1 + J_2 + J_3, \\
h(J_3) = J_3 - J_2, \\
h(J_4) = \sum\limits_{i=2}^{k-1} J_i, \\
h(J_{5+i}) = \vp_{\tha_{k+n}}^i([k, \tha_{k+n}(k)])$, $0 \le i \le k-6, \\
h(J_j) = \vp_{\tha_{k+n}}^{k-5}(J_j)$, $k \le j \le k+n-1. 
\end{cases}
$$

\noindent
Then it is easy to see that $(h \circ \vp_{\al_{k+n}})(J_i) = (\vp_{\tha_{k+n}} \circ h)(J_i)$ for every integer $1 \le i \le k+n-1$.  Thus, $h \circ \vp_{\al_{k+n}} = \vp_{\tha_{k+n}} \circ h$ on $W_{{{\mathbb{R}}}^{k+n-1}}$.  Since $h$ is an isomorphism from $W_{{{\mathbb{R}}}^{k+n-1}}$ onto itself, this implies that $\vp_{\al_{k+n}}$ is conjugate to $\vp_{\tha_{k+n}}$ through $h$.  Consequently, the Petrie matrix of $\al_{k+n}$ is similar to that of $\tha_{k+n}$.  This shows that $\al_k$ and $\tha_k$ are right similar and (by taking $n = 0$) the Petrie matrices of $\al_5$ and $\tha_5$ are similar.  On the other hand, for $k > 5$, the Petrie matrices of $\al_k$ and $\tha_k$ are not similar because they have distinct traces.  

Let $\mu_{k+2}$ and $\nu_{k+2}$ denote the cyclic permutations $1 \ra 4 \ra 3 \ra 6 \ra 7 \ra 8 \ra \cdots \ra k \ra k+1 \ra k+2 \ra 5 \ra 2 \ra 1$ and $1 \ra k+1 \ra k+2 \ra k \ra k-1 \ra k-2 \cdots \ra 7 \ra 6 \ra 3 \ra 4 \ra 5 \ra 2 \ra 1$ on $P_{k+2}$ respectively.  By using elementary row operations, it is easy to see that the determinant of the Petrie matrix of $\mu_{k+2}$ is $(-1)^{k+1}$ while that of $\nu_{k+2}$ is $(-1)^k$.  Since they are different, this shows that $\al_k$ and $\tha_k$ are not two-sided similar.

Finally, the Petrie matrices of the cyclic permutations $1 \ra 4 \ra 3 \ra 6 \ra 7 \ra 8 \ra \cdots \ra k-1 \ra k \ra k+1 \ra 5 \ra 2 \ra 1$ and $1 \ra k+1 \ra k \ra k-1 \ra \cdots 7 \ra 6 \ra 3 \ra 4 \ra 5 \ra 2 \ra 1$ are not similar because they have distinct traces.  So, $\al_k$ and $\tha_k$ are not left similar.
\hfill\sq
\bigskip

\begin{figure}[htb]
\centerline{\epsfig{file=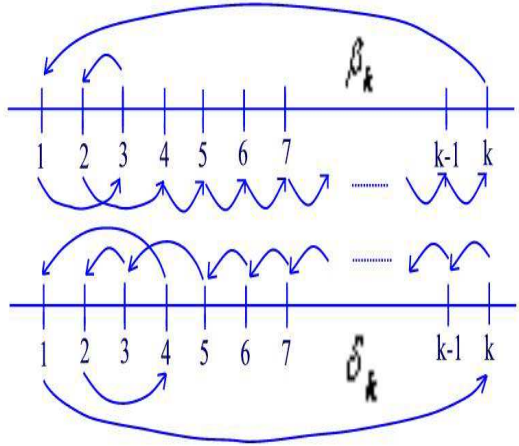,width=7cm,height=6cm}}
\caption{The cyclic permutations $\be_k$ and $\dl_k$ in Theorem 13.}
\end{figure}

\bigskip
\noindent
{\bf Theorem 13.}
{\it Let $k \ge 5$ be an integer.  Then the cyclic permutations \\
\indent\indent\indent $\be_k : 1 \ra 3 \ra 2 \ra 4 \ra 5 \ra 6 \ra 7 \ra \cdots \ra k-1 \ra k \ra 1$ and \\
\indent\indent\indent $\dl_k : 1 \ra k \ra k-1 \ra k-2 \ra \cdots \ra 6 \ra 5 \ra 3 \ra 2 \ra 4 \ra 1$ \\
on $P_k$ are right, but neither left nor two-sided, similar.  Moreover, their respective Petrie matrices are not similar because they have distinct traces.}

\noindent
{\it Proof.}
It is easy to see that the Petrie matrix of $\be_k$ is not similar to that of $\dl_k$ because they have distinct traces.  We now show that $\be_k$ and $\dl_k$ are right similar.  

For every positive integer $n$, let $(\be_{k+n}, \dl_{k+n})$ be a synchronized right extension of $\be_k$ and $\dl_k$ on the set $P_{k+n}$.  Let $\vp_{\be_{k+n}}$ and $\vp_{\dl_{k+n}}$ be the linear transformations on $W_{{{\mathbb{R}}}^{k+n-1}}$ determined by the linearizations of $\be_{k+n}$ and $\dl_{k+n}$ respectively as introduced in section 1.  Let $S = \{ \, J_1, J_1+J_2+J_3$, $\sum\limits_{i=4}^{k-1} J_i$, $[k, \dl_{k+n}(k)]$, $\vp_{\dl_{k+n}}([k, \dl_{k+n}(k)])$, $\vp_{\dl_{k+n}}^2([k, \dl_{k+n}(k)])$, $\vp_{\dl_{k+n}}^3([k, \dl_{k+n}(k)])$, $\cdots, \vp_{\dl_{k+n}}^{k-5}([k, \dl_{k+n}(k)])$, $\vp_{\dl_{k+n}}^{k-4}(J_k)$, $\vp_{\dl_{k+n}}^{k-4}(J_{k+1})$, $\vp_{\dl_{k+n}}^{k-4}(J_{k+2})$, $\cdots, \vp_{\dl_{k+n}}^{k-4}(J_{k+n-1}) \, \}$ and let $V$ be the vector space spanned by $S$.  We want to show that $V = W_{{{\mathbb{R}}}^{k+n-1}}$.  

By direct computations, if $k = 5$, then $\{ \, \vp_{\dl_{5+n}}(J_i) : 1 \le i \le 4 \, \} = \{ \, J_4, J_2+J_3, J_1, \sum\limits_{i=1}^4 J_i +[5, \dl_{5+n}(5)] \, \} \sub V$.  If $k=6$, then $\{ \, \vp_{\dl_{6+n}}^2(J_i) : 1 \le i \le 5 \, \} = \{ \, \sum\limits_{i=1}^5 J_i +[6, \dl_{6+n}(6)], J_1+J_2+J_3, J_4+J_5, J_2+J_3+J_4+J_5, J_1+ \sum\limits_{i=1}^5 J_i +[6, \dl_{6+n}(6)]+\vp_{\dl_{6+n}}([6, \dl_{6+n}(6)]) \, \} \sub V$.  If $k \ge 7$, then, since $\vp_{\dl_{k+n}}(\sum\limits_{i=1}^3 J_i) = \sum\limits_{i=1}^{k-1} J_i$ and $\vp_{\dl_{k+n}}(\sum\limits_{i=4}^{k-1} J_i) = \sum\limits_{i=1}^{k-1} J_i + [k, \dl_{k+n}(k)]$, we easily obtain that $$\vp_{\dl_{k+n}}(J_1) = \sum\limits_{i=4}^{k-1} J_i \quad \in \quad < J_1, J_2+J_3, \sum\limits_{i=4}^{k-1} J_i, [k, \dl_{k+n}(k)] >,$$  $$\vp_{\dl_{k+n}}^2(J_2) = \sum\limits_{i=1}^3 J_i \quad \in \quad < J_1, J_2+J_3, \sum\limits_{i=4}^{k-1} J_i, [k, \dl_{k+n}(k)] >,$$   $$\vp_{\dl_{k+n}}^2(J_3) = \sum\limits_{i=4}^{k-1} J_i \quad \in \quad < J_1, J_2+J_3, \sum\limits_{i=4}^{k-1} J_i, [k, \dl_{k+n}(k)] >,$$  $$\vp_{\dl_{k+n}}^3(J_4) = \sum\limits_{i=1}^3 J_i+\sum\limits_{i=1}^{k-1} J_i +[k, \dl_{k+n}(k)] \quad \in \quad < J_1, J_2+J_3, \sum\limits_{i=4}^{k-1} J_i, [k, \dl_{k+n}(k)] >.$$  Therefore, $\vp_{\dl_{k+n}}^{k-4}(J_i) \in V$ for all $1 \le i \le 4$.  

On the other hand, for $k \ge 7$ and $2 \le i \le k-5$, we have $\vp_{\dl_{k+n}}^{k-5-i}(J_{k-i}) = J_5$ and, since $$\vp_{\dl_{k+n}}^3(J_5)) = J_2  + J_3 + 2\sum\limits_{i=4}^{k-1} J_i \,\, \in \,\, < J_1, J_2+J_3, \sum\limits_{i=4}^{k-1} J_i, [k, \dl_{k+n}(k)] >,$$  $$\vp_{\dl_{k+n}}(J_2+J_3) = J_1+J_2+J_3, \quad \vp_{\dl_{k+n}}(\sum\limits_{i=1}^3 J_i) = \sum\limits_{i=1}^{k-1} J_i,$$ and $$\vp_{\dl_{k+n}}\biggl(\sum\limits_{i=4}^{k-1} J_i \biggr) = \sum\limits_{i=1}^{k-1} J_i + [k, \dl_{k+n}(k)],$$  we have $\vp_{\dl_{k+n}}^{k-4}(J_{k-i}) = \vp_{\dl_{k+n}}^{i-2}(\vp_{\dl_{k+n}}^3(J_5)) = \vp_{\dl_{k+n}}^{i-2}(J_2  + J_3 + 2\sum\limits_{i=4}^{k-1} J_i)$ $$\in \quad < J_1, J_2+J_3, \sum\limits_{i=4}^{k-1} J_i, [k, \dl_{k+n}(k)], \vp_{\dl_{k+n}}([k, \dl_{k+n}(k)]), \cdots, \vp_{\dl_{k+n}}^{k-7}([k, \dl_{k+n}(k)]) >.$$  
So, we inductively obtain that $\vp_{\dl_{k+n}}^{k-4}(J_{k-i}) \in V$ for all $2 \le i \le k-5$.  That is, $\vp_{\dl_{k+n}}^{k-4}(J_j) \in V$ for all $5 \le j \le k-2$.  Finally, it is easy to see that $\vp_{\dl_{k+n}}^{k-4}(J_{k-1}) = \vp_{\dl_{k+n}}^{k-5}(J_{k-2}+J_{k-1}+[k, \dl_{k+n}(k)]) = \vp_{\dl_{k+n}}^2(\vp_{\dl_{k+n}}^{k-7}(J_{k-2}+J_{k-1}+[k, \dl_{k+n}(k)])) = \vp_{\dl_{k+n}}^2(J_5+J_6+J_7+ \cdots + J_{k-1} + [k, \dl_{k+n}(k)] + \vp_{\dl_{k+n}}([k, \dl_{k+n}(k)]) + \cdots + \vp_{\dl_{k+n}}^{k-7}([k, \dl_{k+n}(k)]) = J_1 + \sum\limits_{i=1}^{k-1} J_i + \sum\limits_{i=0}^{k-5} \vp_{\dl_{k+n}}^i([k, \dl_{k+n}(k)]) \in V$.  

Therefore, we have shown that $V \sub W_{{{\mathbb{R}}}^{k+n-1}}$.  Since, by Theorem 6, $\vp_{\dl_{k+n}}$ is an isomorphism on $W_{{{\mathbb{R}}}^{k+n-1}}$, $\{ \, \vp_{\dl_{k+n}}^{k-4}(J_i) : 1 \le i \le k+n-1 \, \}$ is a basis for $W_{{{\mathbb{R}}}^{k+n-1}}$.  Thus, $V = W_{{{\mathbb{R}}}^{k+n-1}}$.  In particular, $S$ is a basis for $W_{{{\mathbb{R}}}^{k+n-1}}$.

Let $h$ the the linear isomorphism from $W_{{{\mathbb{R}}}^{k+n-1}}$ onto itself defined by
$$
\begin{cases}
h(J_1) = J_1, \\
h(J_2) = J_1 + J_2 + J_3, \\
h(J_3) = \sum\limits_{i=4}^{k-1} J_i, \\
h(J_{4+i}) = \vp_{\dl_{k+n}}^i([k, \dl_{k+n}(k)]), \,\, 0 \le i \le k-5, \\
h(J_j) = \vp_{\dl_{k+n}}^{k-4}(J_j), \,\, k \le j \le k+n-1. 
\end{cases}
$$

Now $(h \circ \vp_{\be_{k+n}})(J_1) = h(J_3) = \sum\limits_{i=4}^{k-1} J_i = \vp_{\dl_{k+n}}(J_1) = (\vp_{\dl_{k+n}} \circ h)(J_1)$, $(h \circ \vp_{\be_{k+n}})(J_2) = h(J_2+J_3) = J_1+J_2+J_3 + \sum\limits_{i=4}^{k-1} J_i = \vp_{\be_{k+n}}(J_1+J_2+J_3) = (\vp_{\be_{k+n}} \circ h)(J_2)$, and $(h \circ \vp_{\be_{k+n}})(J_3) = h(J_2+J_3+J_4) = J_2+J_3+J_4 + \sum\limits_{i=4}^{k-1} J_i + [k, \dl_{k+n}(k)] = \vp_{\rho_{k+n}}(\sum\limits_{i=4}^{k-1} J_i) = (\vp_{\dl_{k+n}} \circ h)(J_3)$.  Furthermore, for $0 \le i \le k-5$, $(h \circ \vp_{\be_{k+n}})(J_{4+i}) = h(J_{4+i+1}) = \vp_{\dl_{k+n}}^{i+1}([k, \dl_{k+n}(k)]) = (\vp_{\dl_{k+n}} \circ h)(J_{4+i})$.  

If $k \le j \le k+n-1$ and $\vp_{\be_{k+n}}(J_j) = [s, t] \subset [k, k+n]$, then $\vp_{\be_{k+n}}(J_j) = [\be_{k+n}(j) : \be_{k+n}(j+1)] = [s, t] = [\dl_{k+n}(j) : \dl_{k+n}(j+1)] = \vp_{\dl_{k+n}}(J_j)$.  So, $(h \circ \vp_{\be_{k+n}})(J_j) = h([s, t]) = h(\sum\limits_{i=s}^{t-1} J_i) = \sum\limits_{i=s}^{t-1} h(J_i) = \sum\limits_{i=s}^{t-1} \vp_{\dl_{k+n}}^{k-4}(J_j) = \vp_{\dl_{k+n}}^{k-4}(\sum\limits_{i=s}^{t-1} J_i) = \vp_{\dl_{k+n}}^{k-4}([s, t]) = \vp_{\dl_{k+n}}^{k-4}(\vp_{\dl_{k+n}}(J_j)) = \vp_{\dl_{k+n}}(\vp_{\dl_{k+n}}^{k-4}(J_j)) = \vp_{\dl_{k+n}}(h(J_j)) = (\vp_{\dl_{k+n}} \circ h)(J_j)$, where $[a : b]$ denotes the closed interval with $a$ and $b$ as endpoints.  

If $k \le j \le k+n-1$ and $\vp_{\be_{k+n}}(J_j) = [1, t]$ with $k < t \le k+n$, then $\vp_{\dl_{k+n}}(J_j) = [k-1, t]$.  So, $(\vp_{\dl_{k+n}} \circ h)(J_j) = \vp_{\dl_{k+n}}(h(J_j)) = \vp_{\dl_{k+n}}(\vp_{\dl_{k+n}}^{k-4}(J_j)) = \vp_{\dl_{k+n}}^{k-4}(\vp_{\dl_{k+n}}(J_j)) = \vp_{\dl_{k+n}}^{k-4}([k-1, t]) = \vp_{\dl_{k+n}}^{k-4}(J_{k-1}+\sum\limits_{i=k}^{t-1} J_i) = \vp_{\dl_{k+n}}^{k-4}(J_{k-1}) + \sum\limits_{i=k}^{t-1} \vp_{\dl_{k+n}}^{k-4}(J_i) = \vp_{\dl_{k+n}}^{k-5}(\vp_{\dl_{k+n}}(J_{k-1})) + \sum\limits_{i=k}^{t-1} \vp_{\dl_{k+n}}^{k-4}(J_i) = \vp_{\dl_{k+n}}^{k-5}(J_{k-2}+J_{k-1}+[k, \dl_{k+n}(k)]) + \sum\limits_{i=k}^{t-1} \vp_{\dl_{k+n}}^{k-4}(J_i) = J_1 + (J_1+J_2+J_3) + \sum\limits_{i=4}^{k-1} J_i + \sum\limits_{i=0}^{k-5} \vp_{\dl_{k+n}}^i([k, \dl_{k+n}(k)])$ $+ \sum\limits_{i=k}^{t-1} \vp_{\dl_{k+n}}^{k-4}(J_i) = h(J_1)+h(J_2)+h(J_3) + \cdots + h(J_{t-1}) = h([1, t]) = (h \circ \vp_{\be_{k+n}})(J_j)$.  

Therefore, we have shown that $(h \circ \vp_{\be_{k+n}})(J_i) = (\vp_{\dl_{k+n}} \circ h)(J_i)$ for every integer $1 \le i \le k+n-1$.  Thus, $h \circ \vp_{\be_{k+n}} = \vp_{\dl_{k+n}} \circ h$ on $W_{{{\mathbb{R}}}^{k+n-1}}$.  Since $h$ is an isomorphism from $W_{{{\mathbb{R}}}^{k+n-1}}$ onto itself, this implies that $\vp_{\be_{k+n}}$ is conjugate to $\vp_{\dl_{k+n}}$ through $h$.  Consequently, the Petrie matrix of $\be_{k+n}$ is similar to that of $\dl_{k+n}$.  This shows that $\be_k$ and $\dl_k$ are right similar.

Let $\mu_{k+2}$ and $\nu_{k+2}$ denote the cyclic permutations $1 \ra 4 \ra 3 \ra 5 \ra 6 \ra 7 \ra \cdots \ra k-1 \ra k \ra k+1 \ra k+2 \ra 2 \ra 1$ and $1 \ra k+1 \ra k+2 \ra k \ra k-1 \ra k-2 \ra \cdots \ra 7 \ra 6 \ra 4 \ra 3 \ra 5 \ra 2 \ra 1$ on $P_{k+2}$ respectively.  By using elementary row operations, it is easy to see that the determinant of the Petrie matrix of $\mu_{k+2}$ is $(-1)^{k+1}$ while that of $\nu_{k+2}$ is $(-1)^k$.  Since they are different, this shows that $\be_k$ and $\dl_k$ are not two-sided similar.

Since the trace of the Petrie matrix of the cyclic permutation $1 \ra 4 \ra 3 \ra 5 \ra 6 \ra 7 \ra \cdots \ra k-1 \ra k \ra 2 \ra 1$ is 5 and that of the Petrie matrix of the cyclic permutation $1 \ra k \ra k-1 \ra k-2 \ra \cdots \ra 7 \ra 6 \ra 4 \ra 3 \ra 5 \ra 2 \ra 1$ is 3, $\be_k$ and $\dl_k$ are not left similar.
\hfill\sq
\bigskip

\section{Equivalence classes of one-sided and two-sided similarities in the symmetric group $S_3$}
In the symmetric group $S_3$ of permutations on the set $\{ 1, 2, 3 \}$, no pairs of permutations can be one-sided similar or weakly similar while the pair $\{ (123), (132) \}$ is the only pair which is two-sided similar.  No other pairs are two-sided weakly similar.  

\section{Equivalence classes of one-sided and two-sided similarities in the symmetric group $S_4$}
In the symmetric group $S_4$ of permutations on the set $\{ 1, 2, 3, 4 \}$, the three pairs $\{ (12), (23) \}$, $\{ (123), (132) \}$, and $\{ (1342), (1432) \}$ are the only pairs which are right similar.  

That the pair $\{ (12), (23) \}$ is right similar is depicted in the following Figure 13.

\begin{figure}[htb]
\centerline{\epsfig{file=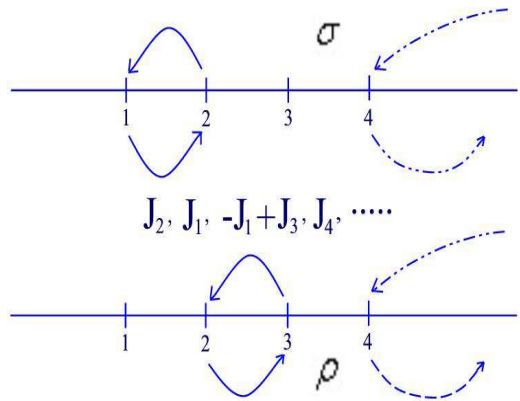,width=7cm,height=6cm}}
\caption{The permutations $\si$ and $\rho$ are right similar.}
\end{figure}

\pagebreak

That the pair $\{ (123), (132) \}$ is right similar is depicted in the following Figure 14.

\begin{figure}[htb]
\centerline{\epsfig{file=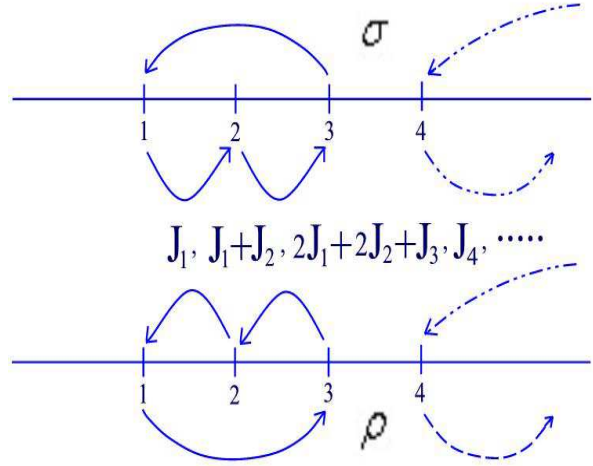,width=8cm,height=7cm}}
\caption{The permutations $\si$ and $\rho$ are right similar.}
\end{figure}

That the pair $\{ (1342), (1432) \}$ is right similar follows from Theorem 7.  

By Theorem 3, the pairs $\{ (1342), (1234) \}$, $\{ (23), (34) \}$, and $\{ (234), (243) \}$ are the only pairs which are left similar.  

As for weakly similarities, the pairs $\{ (34), (12)(34) \}$, $\{ (142), (243) \}$, and the triple $\{ (14)(23)$, $(14), (24) \}$, are the only permutations in $S_4$ which are right weakly similar, but not right similar.  

That the pair $\{ (34), (12)(34) \}$ is not right similar can be seen by considering the synchronized right extension $(\si_5, \rho_5)$ of $\si_4$ and $\rho_4$ on $P_5$ with $\si_5 = (345), \rho_5 = (12)(345)$.  It is easy to check that the minimal polynomial of $M_{\si_5,4}$ is $x^3-2x^2+1$ while that of $M_{\rho_5,4}$ is $x^4-3x^3+2x^2+x-1$.

To see that the pair $\{ (142), (243) \}$ is right weakly similar, let $\si = (142)$ and $\rho = (243)$.  For any $k \ge 5$ and any synchronized right extension $(\si_k, \rho_k)$ of $\si$ and $\rho$, let $M(\si_k,k-1|S)$ ($M(\rho_k,k-1|S)$  respectively) denote the matrix representing the linear transformation $\vp_{\si_k}$ ($\vp_{\rho_k}$ respectively) on $W_{{{\mathbb{R}}}^{k-1}}$ with respect to the basis $S = \{ J_2, J_1, J_2+J_3, J_4, J_5, J_6, \cdots, J_{k-1} \}$.  We then use Laplace's formula to expand the determinants $\det(M(\si_k,k-1|S) - \ld I)$ and $\det(M(\rho_k, k-1 | S) - \ld I)$ along the first column to confirm that they have the same characteristic polynomial.  

To see that the triple $\{ (14)(23)$, $(14), (24) \}$, are right weakly similar to one another, let $\si = (14)(23)$, $\rho = (14)$, and $\gamma = (24)$.  For any $k \ge 5$ and any synchronized right extension $(\si_k, \rho_k, \gamma_k)$ of $\si$, $\rho$, and $\gamma$, let $M(\beta_k,k-1 | S_\beta)$ denote the matrix representing the linear transformation $\vp_{\beta_k}$ on $W_{{{\mathbb{R}}}^{k-1}}$ with respect to the basis $S_\beta$, where 
$$
S_\beta = \begin{cases}
\,\, \{ J_2, J_2+J_3, J_1+J_2+J_3, J_4, J_5, J_6, \cdots, J_{k-1} \}, \,\, \text{if} \,\,\, \beta = \si \,\,\, \text{or} \,\,\, \beta = \rho, \\
\,\, \{ J_1, J_3, J_2+J_3, J_4, J_5, J_6, \cdots, J_{k-1} \}, \,\, \text{if} \,\,\, \beta = \gamma. \\
\end{cases}
$$
We then use Laplace's formula to expand the determinants $\det(M(\si_k,k-1 | S_{\si_k}) - \ld I)$ and $\det(M(\rho_k, k-1 | S_{\rho_k}) - \ld I)$ along the first row and to the determinants $\det(M(\rho_k,k-1 | S_{\rho_k}) - \ld I)$ and $\det(M(\gamma_k, k-1 | S_{\gamma_k}) - \ld I)$ along the first column to confirm that they have the same characteristic polynomial.  

Finally, as for the two-sided similarity or weak similarity, the pair $\{ (134), (142) \}$ and the triple $\{ (1234), (1432), (1342) \}$ are the only permutations in $S_4$ which are two-sided similar and the pair $\{ (23), (14)(23) \}$ is the only pair in $S_4$ which are two-sided weakly similar, but not two-sided similar.  

That the pair $\{ (134), (142) \}$ is two-sided similar can be seen as follows: Let $\si_4 = (134)$ and $\rho_4 = (142)$.  For any synchronized two-sided extension $(\si_{m+4+n}, \rho_{m+4+n})$ of $\si_4$ and $\rho_4$, let $S = \{ J_1, J_2, \cdots, J_m$, $J_{m+1}+J_{m+2}$, $J_{m+3}, [m+4, \rho_{m+4+n}(m+4)]$, $\vp_{\rho_{m+4+n}}(J_{m+4})$, $\vp_{\rho_{m+4+n}}(J_{m+5})$, $\vp_{\rho_{m+4+n}}(J_{m+6})$, $\cdots, \vp_{\rho_{m+4+n}}(J_{m+2+n})$, $\vp_{\rho_{m+4+n}}(J_{m+3+n}) \}$.  By mimicking the proofs as in the previous theorems, we'll obtain the desired result.  That the triple $\{ (1234), (1432), (1342) \}$ are two-sided similar follows from Theorem 8.  

That the pair $\{ (23), (14)(23) \}$ is two-sided weakly similar can be established by expanding the determinants obtained with respect to the standard bases along the row which corresponds to the basis element $J_{m+2}$.

\end{document}